\documentclass[12pt]{amsart}
\usepackage{amsmath,amssymb,amsthm,amscd,amsxtra,eucal,latexsym,mathrsfs}

\newcommand{\nc}{\newcommand}

\newcommand{\cA}{{\mathcal A}}
\newcommand{\cB}{{\mathcal B}}

\newcommand{\cD}{{\mathcal D}}
\newcommand{\cH}{{\mathcal H}}
\newcommand{\cE}{{\mathcal E}}
\newcommand{\cG}{{\mathcal G}}

\newcommand{\cO}{{\mathcal O}}
\newcommand{\cL}{{\mathcal L}}

\newcommand{\cF}{{\mathcal F}}
\newcommand{\cK}{{\mathcal K}}
\newcommand{\cP}{{\mathcal P}}

\newcommand{\cR}{{\mathcal R}}

\newcommand{\cT}{{\mathcal T}}
\newcommand{\cU}{{\mathcal U}}

\newcommand{\cX}{{\mathcal X}}
\newcommand{\cY}{{\mathcal Y}}

\renewcommand{\AA}{{\mathbb A}}

\newcommand{\ZZ}{{\mathbb Z}}
\newcommand{\QQ}{{\mathbb Q}}

\newcommand{\WW}{{\mathbb W}}

\newcommand{\VV}{{\mathbb V}}

\newcommand{\gr}{\mathfrak{r}}
\newcommand{\gs}{\mathfrak{s}}

\newcommand{\on}{\operatorname}

\newcommand{\Rep}{{\on{Rep}}}
\newcommand{\Sch}{{\on{Sch}}}

\newcommand{\Qlb}{\mathbb{\bar Q}_\ell}
\newcommand{\Gm}{\mathbb{G}_m}

\newcommand{\toup}[1]{\stackrel{#1}{\to}}
\newcommand{\hook}[1]{\stackrel{#1}{\hookrightarrow}}

\newcommand{\getsup}[1]{\stackrel{#1}{\gets}}

\newcommand{\Hom}{\on{Hom}}

\nc{\cExt}{\on{{\mathcal E}xt}}

\newcommand{\Ker}{\on{Ker}}
\newcommand{\Coker}{\on{Coker}}

\newcommand{\Aut}{\on{Aut}}

\newcommand{\RG}{\on{R\Gamma}}

\newcommand{\Pic}{\on{Pic}}

\newcommand{\Bun}{\on{Bun}}
\newcommand{\uBun}{\on{\underline{Bun}}}

\newcommand{\Bunt}{\on{\widetilde\Bun}}

\newcommand{\rk}{\on{rk}}
\newcommand{\Spec}{\on{Spec}}

\newcommand{\GL}{\on{GL}}

\newcommand{\pr}{\on{pr}}
\newcommand{\id}{\on{id}}

\newcommand{\QED}{$\square$} 
\newcommand{\Fq}{\mathbb{F}_q}  

\newcommand{\iso}{\xrightarrow{\sim}}

\newcommand{\comp}{\circ}

\renewcommand{\H}{{\on{H}}}   
\newcommand{\R}{\on{R}\!}   

\newcommand{\D}{\on{D}}       

\newcommand{\select}[1]{{\it{#1}}}

\newcommand{\und}[1]{\underline{#1}}

\newcommand{\<}{\langle}
\renewcommand{\>}{\rangle}

\newcommand{\Sph}{\on{Sph}}

\newcommand{\Ind}{\on{Ind}}

\nc{\Perv}{\on{Perv}}
\newcommand{\Div}{\on{Div}}
\newcommand{\cPic}{\on{{\cP}ic}}
\newcommand{\Tors}{\on{{\mathcal{T}}\!ors}}

\newtheorem{Lm}{Lemma}[section]

\newtheorem{Pp}{Proposition}[section]
\newtheorem{Cor}{Corollary}[section]

\newtheorem{Question}{Question}[section]

\theoremstyle{remark}
\newtheorem{Rem}{Remark}[section]

\theoremstyle{definition}
\newtheorem{Def}{Definition}[section]

\newenvironment{Prf}{\par\noindent {\it Proof }}{\QED}
\newcommand{\Step}[1]{\par\noindent{\bf Step {#1}}.}

\begin{document}

\title{Twisted geometric Langlands correspondence for a torus}
\author{Sergey Lysenko}
\address{Institut \'Elie Cartan Nancy (Math\'ematiques) 
 Universit\'e de Lorraine, B.P. 239, F-54506 Vandoeuvre-l\`es-Nancy Cedex, 
 France}
\email{sergey.lysenko@univ-lorraine.fr} 

\begin{abstract} 
Let $T$ be a split torus over local or global function field. The theory of Brylinski-Deligne gives rise to the metaplectic central extensions of $T$ by a finite cyclic group. The representation theory of these metaplectic tori has been developed to some extent in the works of M. Weissman, G. Savin, W. T. Gan, P. McNamara and others. In this paper we propose a geometrization of this theory in the framework of the geometric Langlands program (in the everywhere nonramified case). 
\end{abstract}

\maketitle
\tableofcontents

\section{Introduction}

In this paper we propose a setting for a twisted (nonramified) geometric Langlands correspondence for a split torus in the local and global case. Here `twisted' refers to the quantum Langlands correspondence as outlined in (\cite{St}, \cite{G2}) with the quantum parameter being a root of unity. 

 In \cite{W,W2} M. Weissman has proposed a setting for the representation theory of metaplectic groups over local and global fields. In his approach the metaplectic groups are central extensions of a reductive group by a finite cyclic group coming from Brylinsky-Deligne theory \cite{DB}. Considering only the case of a split torus $T$, we first review the approach by Weissman and prove some related results for automorphic forms over these metaplectic groups in the form we need. This is our subject to geometrize.
 
 In the geometric case we work with $\Qlb$-sheaves to underline a relation with the results on the level of functions (a version for $\cD$-modules should also hold as in \cite{St}). Let $k$ be an algebraically closed field.
In the local case for $F=k((t))$ we consider central extensions $1\to B(\mu_n)\to E\to T(F)\to 1$ coming from Heisenberg $\kappa$-extensions of \cite{BD2}. We construct a category of perverse sheaves equipped with an action of $E$, which is a geometric analog of the corresponding irreducible representation on the level of functions. 

 In the global case we start with a smooth projective curve $X$ over $k$. Let $\Bun_T$ be the stack of $T$-torsors on $X$. We explain how the Brylinski-Deligne data yield a $\mu_n$-gerbe $\Bunt_{T,\lambda}\to\Bun_T$. For an injective character $\zeta\colon \mu_n(k)\to\Qlb^*$ let $\D_{\zeta}(\Bunt_{T,\lambda})$ be the derived category  of $\Qlb$-sheaves on $\Bunt_{T,\lambda}$ on which $\mu_n(k)$ acts by $\zeta$. We define Hecke functors on this category leading to the problem of the corresponding spectral decomposition. Finally, for each spectral parameter $E$ we find a Hecke eigen-sheaf $\cK_E$ corresponding to this parameter and irreducible (over each connected component). It is expected to be unique, but this is not proved yet.
 
  The sheaf $\cK_E$ is a local system. We expect that for $E$ trivial it coincides with the Heisenberg local system constructed by R. Bezrukavnikov, M. Finkelberg and V. Schechtman in \cite{BFS} (in the setting of $\cD$-modules), but we could not check this. Our result should also be related to the equivalence of categories of modules over some sheaves of twisted differential operators on abelian varieties obtained in \cite{PR}. 
    
\section{Main results}
\label{section_Main_results}

\subsubsection{Notations} 
\label{sect_notations} 
For a central extension of groups $1\to A\to E\to G\to 1$ we denote for $a,b\in G$ by $(a,b)_c\in E$ the commutator. If $\tilde a, \tilde b\in E$ are over $a,b\in G$ then
$$
(a,b)_c=\tilde a\tilde b\tilde a^{-1}\tilde b^{-1}.
$$ 
Call a map $c\colon G\times G\to A$  \select{alternating} if $c(a,a)=1$ for all $a\in G$.  If $G$ is abelian then $(a,b)_c\in A$, and the map $(\cdot, \cdot)_c\colon G\times G\to A$ is bimultiplicative, alternating and satisfies $(a,b)_c(b,a)_c=1$ for $a,b\in G$.   

 Let $k$ be an algebraically closed field of characteristic $p\ge 0$ (everywhere except in Section~\ref{sec_classical}, where we assume $k=\Fq$). Pick a prime $\ell$ invertible in $k$. Let $\Qlb$ denote an algebraic closure of $\QQ_{\ell}$. 
All our schemes or stacks are defined over $k$. For an algebraic stack $S$ locally of finite type write $\D(S)$ for the category introduced in (\cite{LO}, Remark~3.21) and denoted $\D_c(S,\Qlb)$ in \select{loc.cit.} It should be thought of as the unbounded derived category of constructible $\Qlb$-sheaves on $S$. Our convention is that a super line is a $\ZZ/2\ZZ$-graded line.
 
 Let $X$ be a smooth projective absolutely irreducible curve over $k$. Write $g$ for the genus of $X$. For an algebraic group $G$ over $k$ write $B(G)$ for the stack quotient of $\Spec k$ by $G$. Denote by $\Bun_G$ the stack classifying $G$-torsors on $X$. For a split torus $T$, a $T$-torsor $\cF$ on a base $S$ and a character $\check{\lambda}$ of $T$ denote by $\cL^{\check{\lambda}}_{\cF}$ the line bundle on $S$ corresponding to the push-forward of $\cF$ via $\check{\lambda}: T\to\Gm$. 
 
  For a notion of a group stack and an action of a group stack on another stack we refer the reader to (\cite{GL}, Appendix A). 
  
\subsubsection{} Let $\Lambda$ be a free abelian group of finite type, set $\check{\Lambda}=\Hom(\Lambda,\ZZ)$ and $T=\Lambda\otimes\Gm$.

In Section~3 we assume $k$ finite, let $\AA$ be the adeles ring of $k(X)$. We consider a central extension $\bar E$ of $T(\AA)$ by a finite cyclic group as in \cite{W} coming from the Brylinski-Deligne theory \cite{DB}. We present basic results about the representation theory of $\bar E$ in the everywhere nonramified case. They are partially borrowed and partially inspired by the works of M. Weissman \cite{W, W2}. This is our subject to geometrize. 

 Starting from Section~\ref{section_Preliminaries} we assume $k$ algebraically closed. Recall that $\Bun_T$ denotes the stack of $T$-torsors on $X$. This is a commutative group stack. In Section~\ref{section_Preliminaries} we review a relation between gerbes on a given base $Z$ and central extensions of the fundamental group of $Z$. We also review and introduce notations for the category of $\theta$-data of Beilinson-Drinfeld \cite{BD2} and review its relation to the central extensions of $\Bun_T$ by $\Gm$. 
 
\subsubsection{} In Section~\ref{section_Geomerization} we consider the twisted geometric Langlands correspondence for a torus in the geometric setting. In the local case we let $\cO=k[[t]]\subset F=k((t))$. Let $\kappa\colon \Lambda\otimes\Lambda\to\ZZ$ be an even symmetric bilinear form. Pick $n>1$ invertible in $k$, let $\zeta\colon\mu_n(k)\to\Qlb^*$ be an injective character. The Heisenberg $\kappa$-extensions of $T(F)$ (\cite{BD}, Definition~10.3.13) give rise to a central extensions of group stacks
\begin{equation}
\label{ext_intro_T(F)_by_B(mu_n)}
1\to B(\mu_n)\to E\to T(F)\to 1.
\end{equation}
We introduce a category of some $\Qlb$-perverse sheaves equipped with an action of $E$, which should play the role of a unique irreducible representation of $E$ with a given central character extending $\zeta$. The central extension (\ref{ext_intro_T(F)_by_B(mu_n)}) splits over $T(\cO)$, we also describe the geometric analog of the corresponding nonramified Hecke algebra (as in \cite{R}). 
 
 In the global case our input data are as follows. Let $\theta=(\kappa,\lambda,c)$ be an object of the category $\cP^{\theta}(X,\Lambda)$ of theta-data (cf. Section~\ref{section_421} and \ref{section_Global setting}). Assume $\kappa\colon\Lambda\otimes\Lambda\to\ZZ$ even. It gives rise to a line bundle also denoted $\lambda$ on $\Bun_T$, it is purely of parity zero (cf. Section~\ref{section_Global setting}). We assume that the associated map $\delta_{\lambda}\colon \Lambda\to \check{\Lambda}$ (defined in Section~\ref{section_421})  equals $\kappa$. Pick $n>1$ and an injective character $\zeta\colon \mu_n(k)\to\Qlb^*$. 

Let $\Bunt_{T,\lambda}$ be the gerbe of $n$-th roots of $\lambda$ over $\Bun_T$. Let 
$$
\Lambda^{\sharp}=\{\mu\in\Lambda\mid \kappa(\mu,\nu)\in n\ZZ\;\mbox{for all}\; \nu\in\Lambda\}.
$$ 
Set $T^{\sharp}=\Lambda^{\sharp}\otimes \Gm$, let $i\colon T^{\sharp}\to T$ be the corresponding isogeny. Let $i_X\colon \Bun_{T^{\sharp}}\to\Bun_T$ be the corresponding push forward map, and $\Bunt_{T^{\sharp},\lambda}$ the restriction of the gerbe $\Bunt_{T, \lambda}$ under $i_X$. Let $\kappa^{\sharp}$ (resp., $\theta^{\sharp}\in\cP^{\theta}(X,\Lambda^{\sharp})$) denote the restriction of $\kappa$ (resp., of $\theta$) to $\Lambda^{\sharp}$. 

 As a part of our input data, we pick a $n$-th root of $\theta^{\sharp}$ in $\cP^{\theta}(X,\Lambda^{\sharp})$. It always exists under our assumptions and gives rise to a section $\gs\colon \Bun_{T^{\sharp}}\to \Bunt_{T^{\sharp},\lambda}$ of the gerbe $\Bunt_{T^{\sharp},\lambda}\to\Bun_{T^{\sharp}}$. For $\mu\in\Lambda$ we denote by $\Bun_T^{\mu}$ the connected component of $\Bun_T$ defined in Section~\ref{section_421}. Let $\Bunt_{T,\lambda}^{\mu}$ be the preimage of $\Bun_T^{\mu}$ in $\Bunt_{T,\lambda}$.
 
  Let $\D_{\zeta}(\Bunt_{T,\lambda}^{\mu})$ be the bounded derived category of $\Qlb$-sheaves on $\Bunt_{T,\lambda}^{\mu}$, on which $\mu_n(k)$ acts via $\zeta$. We have used here the natural action of $\mu_n(k)$ on $\Bunt_{T,\lambda}$ by 2-automorphisms. Write $\D_{\zeta}(\Bunt_{T,\lambda})$ for the derived category of objects whose restriction to each connected component $\Bunt_{T,\lambda}^{\mu}$ lies in $\D_{\zeta}(\Bunt_{T,\lambda}^{\mu})$. For $\mu\in\Lambda^{\sharp}$ define
$\D_{\zeta}(\Bunt_{T^{\sharp},\lambda}^{\mu})$ similarly, we also get $\D_{\zeta}(\Bunt_{T^{\sharp},\lambda})$ as above. 

 The map $\gs$ yields an equivalence $\gs^*\colon \D_{\zeta}(\Bunt_{T^{\sharp},\lambda})\iso \D(\Bun_{T^{\sharp}})$. Our main results are Propositions~\ref{Pp_intro_vanishing}, \ref{Pp_H_eigensheaves} below, their proof is found in Sections~\ref{Section_522}-\ref{Section_524}.  
 
\begin{Pp} 
\label{Pp_intro_vanishing}
Let $\mu\in\Lambda$ with $\mu\notin \Lambda^{\sharp}$. Then $\D_{\zeta}(\Bunt_{T,\lambda}^{\mu})$ vanishes.
\end{Pp}
   
 We define an action on $\Bun_{T^{\sharp}}$ on $\Bunt_{T,\lambda}$ in Section~\ref{section_Hecke_functors}. Let $\check{T}^{\sharp}$ be the torus dual to $T^{\sharp}$ over $\Qlb$. Let $E$ be a $\check{T}^{\sharp}$-local system on $X$. We define a notion of a $E$-Hecke eigen-sheaf in $\D_{\zeta}(\Bunt_{T,\lambda})$. The geometric Langlands problem in our setting is to find a spectral decomposition of $\D_{\zeta}(\Bunt_{T,\lambda})$ under the action of $\Bun_{T^{\sharp}}$. 

  The map $i_X$ lifts to a morphism $\pi\colon \Bunt_{T^{\sharp},\lambda}\to \Bunt_{T,\lambda}$ defined in Section~\ref{Section_522}. Let $AE$ be the automorphic local system on $\Bun_{T^{\sharp}}$ associated to $E$ (cf. Section~\ref{section_Hecke_functors}). Let $W\in \D_{\zeta}(\Bunt_{T^{\sharp},\lambda})$ be an object equipped with $\gs^*W\iso AE$. 
\begin{Pp} 
\label{Pp_H_eigensheaves}
(i) The local system $\pi_!W\in \D_{\zeta}(\Bunt_{T,\lambda})$ is equipped with a structure of a $E$-Hecke eigen-sheaf. \\
(ii) There is $\cK_E\in \D_{\zeta}(\Bunt_{T,\lambda})$, which is an irreducible local system over $\Bunt_{T,\lambda}^{\mu}$ for $\mu\in\Lambda^{\sharp}$ and vanishing over $\Bunt_{T,\lambda}^{\mu}$ for $\mu\notin \Lambda^{\sharp}$ with the following property. There is a $\Qlb$-vector space $\VV$ such that 
$\pi_!W\iso \cK_E\otimes \VV$ over each connected component $\Bunt_{T,\lambda}^{\mu}$, $\mu\in\Lambda^{\sharp}$. Moreover, $\cK_E$ admits a structure of a $E$-Hecke eigen-sheaf in $\D_{\zeta}(\Bunt_{T,\lambda})$. 
\end{Pp}  

\begin{Rem} As in \cite{W}, pick two bases $(\epsilon_i), (\eta_i)$ of $\Lambda$ with $1\le i\le r$ such that $\kappa(\epsilon_i, \eta_j)$ vanishes unless $i=j$, and $\kappa(\epsilon_i, \eta_i)=d_i$, where $d_i$ divides $d_{i+1}$ for all $i$ (as soon as both $d_i$ and $d_{i+1}$ are not zero). Let $e_i$ be the smallest positive integer such that $d_i e_i\in n\ZZ$. Then $\Lambda^{\sharp}=\bigoplus_{i=1}^r (e_i\ZZ)\epsilon_i$. Let $e=\prod_i e_i$ be the order of $\Lambda/\Lambda^{\sharp}$. The rank of $\cK$ is $e^{g}$ and $\dim(\VV)=e^g$.
\end{Rem}

It is natural to ask the following.
\begin{Question}
\label{Con_one}
Is there an equivalence of categories $\D(\Bun_{T^{\sharp}})\iso \D_{\zeta}(\Bunt_{T,\lambda})$ commuting with the action of Hecke functors? For each $\mu\in\Lambda^{\sharp}$ it would identify $\D(\Bun_{T^{\sharp}}^{\mu})$ with $\D_{\zeta}(\Bunt_{T,\lambda}^{\mu})$.
\end{Question}  

According to the quantum Langlands conjectures \cite{St}, in the setting of $\cD$-modules the expected answer to the same question is `yes'. However, the corresponding equivalence is expected to be given by a $\cD$-module on $\Bun_{T^{\sharp}}\times\Bunt_{T,\lambda}$, which is not holonomic, and in 
our $\ell$-adic setting we expect the answer to be negative.

On one hand, we show in Proposition~\ref{Pp_scalar_products} that the `corrected scalar products' of automorphic sheaves are the same in both categories (this is an analog of the Rankin-Selberg convolutions from \cite{L} in the setting of metaplectic tori). On the other hand, in Section~\ref{Section_528}
we give an argument in the case of an elliptic curve and $T=\Gm$ showing that these categories are not expected to be equivalent.  
 
\section{Representations of the metaplectic tori on the level of functions}
\label{sec_classical}

\subsection{Central extensions} 
\label{section_301}
Throughout Section~\ref{sec_classical} we assume $k=\Fq$.
Let $A,G$ be abelian groups, write $A$ multiplicatively and $G$ additively. A normalized 2-cocycle is a map $f\colon G\times G\to A$ satisfying $1=f(0, g)=f(g,0)$ for $g\in G$, and 
$$
f(g_1, g_2)f(g_1+g_2, g_3)=f(g_2, g_3)f(g_1, g_2+g_3)
$$ 
for $g_i\in G$. Such a cocycle gives rise to a central extension $1\to A\to A\times G\to G\to 1$ with the group law
$$
(a_1, g_1)(a_2, g_2)=(a_1a_2f(g_1,g_2), g_1+g_2),  \;\;\;\; a_i\in A, g_i\in G.
$$ 
If we need to emphasize the dependence on $f$, we denote this group by $(A\times G)_f$. For $g_i\in G$ the commutator in $(A\times G)_f$ is given by
\begin{equation}
\label{cocycle_commutator_relation}
(g_1, g_2)_c=\frac{f(g_1, g_2)}{f(g_2, g_1)}\ \ .
\end{equation}  
Note that if $f\colon G\times G\to A$ is any bilinear map then $f$ is a normalized 2-cocycle.

 Let $T^0$ be the group of all maps of sets $h\colon G\to A$ such that $h(0)=1$, let $T^1$ be the group of normalized 2-cocycles $f\colon G\times G\to A$. We have a complex $T^0\to T^1$, where $h$ is mapped to $f$ given by 
$$
f(g_1,g_2)=\frac{h(g_1+g_2)}{h(g_1)h(g_2)} \ \  .
$$  
Let $f,f'\in T^1$. Then given $h\in T^0$ with
$$
\frac{f'(g_1, g_2)}{f(g_1,g_2)}=\frac{h(g_1+g_2)}{h(g_1)h(g_2)}\;\;\;\mbox{for all}\;\; g_i\in G,
$$ 
we get an isomorphism of central extensions $(A\times G)_f\iso (A\times G)_{f'}$ given by  $(a,g)\mapsto (ah(g), g)$. 

 The Picard category of central extensions of $G$ by $A$ is canonically equivalent to the Picard category associated to the complex $T^0\to T^1$ (see \cite{DB}). 
 
\subsection{Brylinski-Deligne extensions} 

\subsubsection{}
\label{section_311}
Let $\Lambda$ be a free abelian group of finite type. Recall that isomorphism classes of central extensions $1\to A\to E\to\Lambda\to 1$ are in bijection with the set of bilinear alternating maps $\Lambda\times\Lambda\to A$, the map associated to a central extension being its commutator. Set $\check{\Lambda}=\Hom(\Lambda,\ZZ)$. 

 Let $T=\Lambda\otimes\Gm$. View $T$ as a scheme over $k=\Fq$. Let $\Sch/k$ be the category of $k$-schemes of finite type equipped with the Zariski topology. The $n$-th Quillen's $K$-theory group of a scheme form a presheaf on $\Sch/k$ as the scheme varies. We denote by $K_2$ the associated sheaf on $\Sch/k$ for the Zariski topology. Let $\cExt(T, K_2)$ denote the category of central extensions of $T$ by $K_2$ (in the category of sheaves of groups on $\Sch/k$). 
 
 A bilinear form $\kappa\colon \Lambda\otimes\Lambda\to\ZZ$ is even if $\kappa(x,x)\in 2\ZZ$ for all $x\in\Lambda$. Recall that quadratic forms $q\colon\Lambda\to\ZZ$ are in bijection with the set of even symmetric bilinear forms $\kappa\colon \Lambda\otimes\Lambda\to\ZZ$. The form associated to $q$ is $\kappa$, where $\kappa(x_1,x_2)=q(x_1+x_2)-q(x_1)-q(x_2)$. 
  
\begin{Pp}[\cite{DB}, Theorem~3.16] There is an equivalence of Picard groupoids between $\cExt(T,K_2)$ and the groupoid $\cE(T, K_2)$ of pairs: an even symmetric bilinear form $\kappa\colon {\Lambda\otimes\Lambda}\to\ZZ$ and a central extension
$1\to k^*\to  \tilde\Lambda\to \Lambda\to 1$ whose commutator is given by $(\mu_1,\mu_2)_c=(-1)^{\kappa(\mu_1,\mu_2)}$, $\mu_i\in\Lambda$.
\end{Pp}

 Let us now discuss some realizations of the objects of $\cE(T,K_2)$. It is convenient for this to introduce a category $\cE^s(T)$ of pairs: a symmetric bilinear form $\kappa\colon\Lambda\otimes\Lambda\to\ZZ$ and a central super extension 
\begin{equation}
\label{ext_tilde_Lambda^s} 
 1\to k^*\to \tilde\Lambda^s\to \Lambda\to 1
\end{equation} 
(as in \cite{BD}, Lemma~3.10.3.1) such that its commutator is $(\gamma_1,\gamma_2)_c=(-1)^{\kappa(\gamma_1,\gamma_2)}$. This means that for every $\gamma\in\Lambda$ we are given a $\ZZ/2\ZZ$-graded (or super) line $\epsilon^{\gamma}$, and for every $\gamma_1,\gamma_2\in\Lambda$ a $\ZZ/2\ZZ$-graded isomorphism 
$$
c^{\gamma_1,\gamma_2}\colon \epsilon^{\gamma_1}\otimes\epsilon^{\gamma_2}\iso \epsilon^{\gamma_1+\gamma_2}
$$ 
such that $c$ is associative, i.e., 
$$
c^{\gamma_1, \gamma_2+\gamma_3}(\id_{\epsilon^{\gamma_1}}\otimes c^{\gamma_2,\gamma_3})=c^{\gamma_1+\gamma_2, \gamma_3}(c^{\gamma_1,\gamma_2}\otimes\id_{\epsilon^{\gamma_3}}),
$$
and one has $c^{\gamma_1, \gamma_2}=(-1)^{\kappa(\gamma_1,\gamma_2)}c^{\gamma_2,\gamma_1}\sigma$, where $\sigma\colon \epsilon^{\gamma_1}\otimes\epsilon^{\gamma_2}\iso \epsilon^{\gamma_2}\otimes\epsilon^{\gamma_1}$ is the super commutativity constraint. 

 The category $\cE^s(T)$ is a Picard groupoid with respect to the tensor product of central extensions. For an object of $\cE^s(T)$ as above for each $\gamma\in\Lambda$ the parity of $\epsilon^{\gamma}$ is $\kappa(\gamma,\gamma)\!\!\mod 2$. So, $\cE(T,K_2)$ is the full Picard subgroupoid of $\cE^s(T)$ consisting of pairs $(\kappa, \tilde\Lambda^s)$ such that $\kappa$ is even. 

 According to (\cite{BD}, Lemma~3.10.3.1), let $v\colon\Lambda\to \ZZ/2\ZZ$ be a morphism and $B\colon \Lambda\times\Lambda\to\ZZ/2\ZZ$ a bilinear form such that 
$$
\kappa(\gamma_1,\gamma_2)\!\!\!\!\mod 2=v(\gamma_1)v(\gamma_2)+B(\gamma_1,\gamma_2)+B(\gamma_2,\gamma_1)
$$ 
for all $\gamma_i\in\Lambda$. Then taking $\epsilon^{\gamma}=k$ of parity $v(\gamma)$, we get a super central extension (\ref{ext_tilde_Lambda^s}) given by the cocycle $(\gamma_1,\gamma_2)\mapsto (-1)^{B(\gamma_1,\gamma_2)}$. 
We will use the following two particular cases of this construction.
   
\begin{Lm} 
\label{Lm_1_realizations}
(1) Let $B\colon\Lambda\otimes\Lambda\to\ZZ$ be a bilinear form, set $\kappa=B+{^tB}$, where $^tB(\gamma_1,\gamma_2)=B(\gamma_2,\gamma_1)$ for $\gamma_i\in\Lambda$. Then $B$ gives rise to the object $(\kappa, \tilde\Lambda^s)\in \cE(T,K_2)$, where $\epsilon^{\gamma}=k$ of parity zero for each $\gamma\in\Lambda$, and the central extension (\ref{ext_tilde_Lambda^s}) is given by the cocycle 
$$
(\gamma_1,\gamma_2)\mapsto (-1)^{B(\gamma_1,\gamma_2)}\  .
$$
(2) Let $\check{\lambda}\in\check{\Lambda}$ and $\kappa=\check{\lambda}\otimes\check{\lambda}$. Let $v\colon \Lambda\to\ZZ/2\ZZ$ be the map $\check{\lambda}\!\!\mod 2$. Consider the super extension $1\to k^*\to \tilde\Lambda^{\check{\lambda}}\to \Lambda\to 1$ given by the zero 2-cocycle, where $\epsilon^{\gamma}=k$ is of parity $v(\gamma)$. Its commutator is 
\begin{equation}
\label{commutator_given_by_kappa}
(\gamma_1,\gamma_2)_c=(-1)^{\kappa(\gamma_1, \gamma_2)}\  .
\end{equation}
(2') Let $\kappa\colon \Lambda\otimes\Lambda\to\ZZ$ be an even symmetric bilinear form. Pick a presentation $\kappa=\sum_i b_i(\check{\lambda}_i\otimes\check{\lambda}_i)$ for some $\check{\lambda}_i\in\check{\Lambda}$. For each $i$ we get a super central extension $\tilde\Lambda^{\check{\lambda}_i}$. The tensor product of their $b_i$-th powers is a central extension $1\to k^*\to \tilde\Lambda\to\Lambda\to 1$ of parity zero with the commutator (\ref{commutator_given_by_kappa}). This is an object of $\cE(T, K_2)$. \QED
\end{Lm}

\subsubsection{} 
\label{section_DB_explcit_cocycles}
In this subsection we describe the Brylinski-Deligne extensions of $T$ by $K_2$ by some explicit cocycles.

Let $B\colon\Lambda\otimes\Lambda\to\ZZ$ be any bilinear form and $\kappa=B+{^tB}$ as in Lemma~\ref{Lm_1_realizations}(1). Recall that $K_1=\cO^*$ as a sheaf on $\Sch/k$ (see \cite{DB}). There is a unique bimultiplicative map $\bar f\colon T\times T\to K_2$ satisfying
$$
\bar f(\lambda_1\otimes c_1, \lambda_2\otimes c_2)=\{c_1, c_2\}^{B(\lambda_1,\lambda_2)}
$$ 
for $\lambda_i\in\Lambda$, $c_i\in\Gm$. Here $\{\cdot,\cdot\}\colon K_1\times K_1\to K_2$ is the product in the graded sheaf of rings $\bigoplus_i K_i$ on $\Sch/k$. It is bilinear and skew-symmetric. The map $\bar f$ is a 2-cocycle, it gives rise to a central extension
\begin{equation}
\label{ext_T_by_K_2}  
1\to K_2\to (K_2\times T)_B\to T\to 1,
\end{equation} 
where $(K_2\times T)_B=K_2\times T$ is equipped with the product $(z_1,t_1)(z_2,t_2)=(z_1z_2\bar f(t_1,t_2), t_1t_2)$, $t_i\in T$, $z_i\in K_2$. We write $\bar f_B=\bar f$ if we need to express the dependence on $B$ of the above cocycle.
 
 By (\cite{DB}, Corollary~3.7), the automorphisms of the extension (\ref{ext_T_by_K_2}) are $\Hom(\Lambda, k^*)$. Namely, $q\in\Hom(\Lambda, k^*)$
defines a homomorphism of sheaves $\bar q\colon T\to K_2$ on $\Sch/k$ such that $\bar q(\lambda\otimes c)=\{c, q(\lambda)\}$ for $c\in\Gm$, $\lambda\in\Lambda$. Note that $\bar q\in \H^0(T, K_2)$.
 
  As in \cite{DB}, define a map $\check{\Lambda}\otimes\check{\Lambda}\to \H^0(T, K_2)$ as follows. Given $\check{\lambda}_i\in \check{\Lambda}$, view $\check{\lambda}_i\in \H^0(T, \cO^*)=\H^0(T, K^1)$, so $\{\check{\lambda}_1, \check{\lambda}_2\}\in \H^0(T, K_2)$ via the product in $\bigoplus_j K_j$. Since the product $K_1\times K_1\to K_2$ is bimultiplicative, the map $\check{\lambda}_1\otimes\check{\lambda}_2\mapsto \{\check{\lambda}_1, \check{\lambda}_2\}$
extends to a linear map $\check{\Lambda}\otimes\check{\Lambda}\to \H^0(T, K_2)$ that we denote $B\mapsto \bar B$. 

 Note that if $\{\epsilon_i\}$, $1\le i\le r$ is a basis of $\Lambda$, $c_i\in \Gm$ then 
$$
\bar B(\prod_{i=1}^r (\epsilon_i\otimes c_i))=\prod_{1\le i,j\le r} \{c_i, c_j\}^{B(\epsilon_i, \epsilon_j)}\  .
$$  
The product in $K_2$ is written multiplicatively. Write the operation in the abelian group $\H^0(T, K_2)$ additively, so $q$ and $B$ define $\bar B+\bar q\in \H^0(T, K_2)$.

 To any $s\in \H^0(T, K_2)$ one associates a 2-cocyle $ds\colon \H^0(T\times T, K_2)$ by 
$$
(ds)(t_1, t_2)=\frac{s(t_1t_2)}{s(t_1)s(t_2)}
$$ 
for $t_i\in T$. For $\bar B+\bar q$ as above we get $d(\bar B+\bar q)=d(\bar B)$.   

\begin{Lm} For a bilinear form $C\colon \Lambda\otimes\Lambda\to \ZZ$ and $B=C-{^tC}$, one has $d(\bar C)=\bar f_B$.
\end{Lm}
\begin{Prf}
Pick a base $\{\epsilon_i\}$ in $\Lambda$. Let $t_1=\prod_i (\epsilon_i\otimes g_i)$, $t_2=\prod_i (\epsilon_i\otimes g'_i)$ with $\lambda_i\in\Lambda$, $g_i, g'_i\in\Gm$. We get
$$
\bar C(t_1t_2)=\bar C(\prod_i (\epsilon_i\otimes g_ig'_i))=\prod_{i,j} \{g_ig'_i, g_jg'_j\}^{C(\epsilon_i, \epsilon_j)}
$$
and
$$
\bar C(t_1)=\prod_{i,j} \{g_i, g_j\}^{C(\epsilon_i, \epsilon_j)},
$$
$$
\bar C(t_2)=\prod_{i,j} \{g'_i, g'_j\}^{C(\epsilon_i, \epsilon_j)}\ .
$$
So, 
$$
\frac{\bar C(t_1t_2)}{\bar C(t_1)\bar C(t_2)}=\prod_{i,j} \{g_i, g'_j\}^{C(\epsilon_i, \epsilon_j)-C(\epsilon_j, \epsilon_i)}\ .
$$
Our claim follows.
\end{Prf}

\medskip

 We get a complex $\check{\Lambda}\otimes\check{\Lambda}\toup{d_1} \check{\Lambda}\otimes\check{\Lambda}\toup{d_2} Q(\Lambda)$ say in degrees 0,1,2, where $Q(\Lambda)$ is the group of even symmetric bilinear forms $\Lambda\otimes\Lambda\to \ZZ$. The map $d_1$ is given by $d_1(C)=C-{^tC}$ and $d_2(B)=B+{^tB}$. This complex is exact in degree 1.
 
By (\cite{DB}, Theorem~3.16), the isomorphism class of the extension $(K_2\times T)_B$ is given by $d_2(B)\in Q(\Lambda)$. If $d_2(B)=0$, an element $C\in \check{\Lambda}\otimes\check{\Lambda}$ with $d_1(C)=B$ gives a trivialization of $(K_2\times T)_B$.

\subsection{Global setting}
 
\subsubsection{} 
\label{section_321}
Let $X$ be a smooth projective absolutely irreducible curve over $k$. Set $F=k(X)$, let $\AA$ denote the adeles ring of $X$, $\cO\subset\AA$ be the integer adeles. For $x\in X$ write $F_x$ for the completion of $F$ at $x$. Let $\cO_x\subset F_x$ be the ring of integers and $k(x)$ its residue field. 
 
Write $v(f)$ for the valuation of $f\in F_x^*$. For $f,g\in F_x^*$ write $(\cdot, \cdot)_{st}$ for the tame symbol given by 
$$ 
(f,g)_{st}=(-1)^{deg(x)v(f)v(g)}N_{k(x)/k}((g^{v(f)}f^{-v(g)})(x))
$$
according to (\cite{MP_PR}, Remark~2.2). Here $N_{k(x)/k}$ denotes the norm for the extension $k\subset k(x)$. The global tame symbol $(\cdot, \cdot)_{st}\colon \AA^*\times\AA^*\to k^*$ is given by
$$
(a,b)_{st}=\prod_{x\in X} (a_x, b_x)_{st}\  .
$$
If both $a_x,b_x\in \cO_x$ then $(a_x, b_x)_{st}=1$, so the product is finite. 
 
 Let $B\colon \Lambda\otimes\Lambda\to\ZZ$ be any bilinear form. Set $\kappa=B+{^tB}$, so $\kappa$ is even and symmetric. As in Section~\ref{section_DB_explcit_cocycles}, there is a unique bimultiplicative map $f\colon T(\AA)\times T(\AA)\to k^*$ such that for $\lambda_i\in\Lambda$, $g_i\in \AA^*$ one has
$$
f(\lambda_1\otimes g_1, \lambda_2\otimes g_2)=(g_1,g_2)_{st}^{B(\lambda_1,\lambda_2)}\ .
$$
The map $f$ is a normalized 2-cocycle as defined in Section~\ref{section_301}. 
Let $E=k^*\times T(\AA)$ with the product $(a_1, g_1)(a_2, g_2)=(a_1a_2f(g_1, g_2), g_1g_2)$ for $a_i\in k^*, g_i\in T(\AA)$. Then $E$ is a locally compact topological group that fits into a central extension
\begin{equation}
\label{ext_T(AA)_by_k^*}
1\to k^*\to E\to T(\AA)\to 1.
\end{equation}
The commutator in $E$ is given by 
$$
(\lambda_1\otimes g_1, \lambda_2\otimes g_2)_c=(g_1, g_2)_{st}^{\kappa(\lambda_1, \lambda_2)}
$$
for $g_i\in \AA^*, \lambda_i\in\Lambda$. 

 The map $T(\cO)\to E$, $z\mapsto (1,z)$ is a group homomorphism, a splitting of (\ref{ext_T(AA)_by_k^*}) over $T(\cO)$. The image of $T(\cO)$ is usually not a normal subgroup in $E$. In particular, if $\kappa$ is non-degenerate then $T(\cO)$ is not a normal subgroup.
 
 If $g_i\in F^*$ then $(g_1, g_2)_{st}=1$ (see \cite{ACK}). So, the map $T(F)\to E$, $z\mapsto (1,z)$ is a group homomorphism, a splitting of (\ref{ext_T(AA)_by_k^*}) over $T(F)$.  The two splittings coincide over $T(F)\cap T(\cO)=T(k)$. 
If $\kappa$ is non-degenerate then $T(F)$ is not normal in $E$. 
 
 Let $\Bun_T$ be the stack of $T$-torsors on $X$. Recall that $\Bun_T(k)\iso T(F)\backslash T(\AA)/T(\cO)$ naturally, so $T(F)\backslash E/T(\cO)\to\Bun_T(k)$ is a $k^*$-torsor. The problem of the geometrization of this torsor is essentially solved in (\cite{BD2}, Proposition~3.10.7.1, Lemma~3.10.3.1 and Proposition~4.9.1.2). It will be discussed below in Section~\ref{section_Geomerization}. 
 
\subsubsection{} Pick $n\ge 1$. We assume $n\mid q-1$, so both $\mu_n(k)$ and $k^*/(k^*)^n$ are cyclic of order $n$. The natural map $\mu_n(k)\to k^*/(k^*)^n$ is not always an isomorphism. 

 Denote by
\begin{equation}
\label{ext_EnB_T(AA)_by_mu_n}
1\to  k^*/(k^*)^n\to \bar E\to T(\AA)\to 1
\end{equation} 
the push-forward of (\ref{ext_T(AA)_by_k^*}) by $k^*\to k^*/(k^*)^n$. Let
\begin{equation}
\label{ext_BnEx_T(AA)_by_mu_n}
1\to  k^*/(k^*)^n\to \bar E_x\to T(F_x)\to 1
\end{equation}
be the pull-back of (\ref{ext_EnB_T(AA)_by_mu_n}) by $T(F_x)\to T(\AA)$. 

 The results of this subsection are partially borrowed from and partially inspired by the papers by M. Weissman \cite{W, W2}. Let 
\begin{equation}
\label{def_Lambda^sharp}
\Lambda^{\sharp}=\{\mu\in\Lambda\mid \kappa(\mu,\nu)\in n\ZZ\;\mbox{for all}\; \nu\in\Lambda\}.
\end{equation}
Set $T^{\sharp}=\Lambda^{\sharp}\otimes\Gm$, so we have an isogeny $T^{\sharp}\to T$. 

 Let $Z^{\dag}\subset T(\AA)$ be the subgroup such that the center $Z(\bar E)$ of $\bar E$ is the preimage of $Z^{\dag}$ in $\bar E$. Let $Z^{\dag}_x\subset T(F_x)$ be the subgroup such that the center $Z(\bar E_x)$ of $\bar E_x$ is the preimage of $Z^{\dag}_x$ in $\bar E_x$.
 
\begin{Lm} 
\label{Lm_description_center}
(1) Let $g=(g_x)\in T(\AA)$. Then $g\in Z^{\dag}$ iff $g_x\in Z^{\dag}_x$ for all $x\in X$.\\
(2) $Z^{\dag}_x$ is the image of $T^{\sharp}(F_x)\to T(F_x)$. 
\end{Lm}
\begin{Prf}
(1) Pick a base $\{\epsilon_i\}$ of $\Lambda$. For $g_i\in \AA^*$ the condition
$g=\prod_i (\epsilon_i \otimes g_i)\in Z^{\dag}$ is equivalent to requiring that for any $u\in \AA^*$, $\lambda\in\Lambda$ one has
$$
\prod_i (g_i, u)_{st}^{\kappa(\epsilon_i, \lambda)} \in (k^*)^n \ .
$$
This condition is local over $X$. \\
(2) is (\cite{W}, Proposition~4.1). The property used by Weissman in \select{loc.cit.} is as follows. Let $d\ge 1$. Write $e$ for the smallest positive integer such that $de\in n\ZZ$. Then 
$$
\{g\in F_x^*\mid \mbox{for all}\; h\in F_x^*, (g,h)_{st}^d\in (k^*)^n\}=(F_x^*)^e \ .
$$
The case $d=1$ follows from the non-degeneracy of the tame symbol. Here is the reduction to the case $d=1$. First, if $d$ is prime to $n$ then let $\alpha,\beta\in\ZZ$ with $d\alpha+n\beta=1$. If $g^d=z^n$ for some $z\in F_x^*$ then $g=g^{d\alpha+n\beta}=(z^{\alpha}g^{\beta})^n\in (F_x^*)^n$. Now let $d,n$ be any, set $a=GCD(n,d)$ and $n=ae$, $d=ad'$. Let $z\in F_x^*$ be such that $g^{d'a}=z^{ae}$. Then $g^{d'}z^{-e}\in \mu_a(k)$. The map $\mu_n(k)\to\mu_a(k)$, $u\mapsto u^e$ is surjective. So, $g^{d'}\in (F_x^*)^e$. Since $(d', e)=1$ we get $g\in (F_x^*)^e$. 
\end{Prf}
 
\begin{Rem} 
\label{Rem_to_be_checked_Kottwitz}
The Poitou-Tate duality in Galois cohomology implies $\{y\in \AA^*\mid (y,z)_{st}\in (k^*)^n\;\mbox{for all}\; z\in F^*\}=F^*(\AA^*)^n$. 
\end{Rem} 

  Recall that (\ref{ext_EnB_T(AA)_by_mu_n}) splits canonically over $T(F)$. 
  
\begin{Lm} 
\label{Lm_global_max_abelain_subgroup}
The group $T(F)Z(\bar E)$ is a maximal abelian subgroup of $\bar E$.
\end{Lm}
\begin{Prf}
\Step 1 Pick two bases $(\epsilon_i), (\eta_i)$ of $\Lambda$ with $1\le i\le r$ such that $\kappa(\epsilon_i, \eta_j)$ vanishes unless $i=j$, and $\kappa(\epsilon_i, \eta_i)=d_i$, where $d_i$ divides $d_{i+1}$ for all $i$. Let 
$e_i$ be the smallest positive integer such that $d_i e_i\in n\ZZ$. In this notation we have  
$$
Z^{\dag}=\{\prod_i (\epsilon_i\otimes g_i)\in T(\AA)\mid g_i\in (\AA^*)^{e_i}\;\mbox{for all}\; i\} \ .
$$
Note that $\Lambda^{\sharp}=\bigoplus_{i=1}^r (e_i\ZZ)\epsilon_i$. Now let $g_i\in \AA^*$ and $g=\prod_i(\epsilon_i\otimes g_i)$. If $g$ commutes with $T(F)$ in $\bar E$ then for each $1\le j\le r$ and $u\in F^*$ we get
$$
(\prod_i(\epsilon_i\otimes g_i), \eta_j\otimes u)_c=(g_j, u)_{st}^{d_j}\in (k^*)^n \ .
$$
By Remark~\ref{Rem_to_be_checked_Kottwitz}, this implies $g_j^{d_j}\in F^*(\AA^*)^n$. We are reduced to the following.

\medskip

\Step 2 Let $g\in \AA^*$, $d\ge 1$, let $e$ be the smallest positive integer such that $de\in n\ZZ$. If $g^d\in F^*(\AA^*)^n$ then $g\in F^*(\AA^*)^e$. Let $g^d=f z^n$ with $f\in F^*, z\in\AA^*$. First, consider the case $(d,n)=1$. In this case pick $\alpha,\beta\in\ZZ$ with $d\alpha+n\beta=1$. Then $g=g^{d\alpha+n\beta}=(fz^n)^{\alpha}g^{n\beta}\in F^*(\AA^*)^n$. Now let $n,d$ be arbitrary, set $a=GCD(d,n)$ and $n=ae, d=ad'$. Since $(g^{d'}z^{-e})^a=f$, by the Grunwald-Wang Theorem (\cite{AT}, Theorem~1, Chapter~IX), there is $f_1\in F^*$ with $f=f_1^a$. So, there is $v\in \AA^*$ with $v^e=g^{d'}z^{-e}f_1^{-1}$ and $g^{d'}\in F^*(\AA^*)^e$. Since $(d',e)=1$ from the relatively prime case we get $g\in F^*(\AA^*)^e$.
\end{Prf}


\begin{Rem} 
\label{Rem_about_finite_field}
Let $d,n\ge 1$, let $e$ be the smallest positive integer with $de\in n\ZZ$. If $a\in k^*$ and $a^d\in (k^*)^n$ then $a\in (k^*)^e$. This is proved as in Lemma~\ref{Lm_description_center} (2).
\end{Rem}

\subsubsection{} 
\label{section_323}
Pick an injective character $\zeta\colon k^*/(k^*)^n\to \Qlb^*$. 
The twisted Langlands correspondence at the level of functions for a torus is the study of the representation
\begin{multline*}
\cR=\{f\colon T(F)\backslash{\bar E}\to\Qlb\mid f(yz)=\zeta(y)f(z)\;\mbox{for all}\; y\in k^*/(k^*)^n\\
\mbox{there is an open subgroup}\; \cK\subset T(\cO),\; f(yu)=f(y), \; u\in \cK\}
\end{multline*}
of $\bar E$ by right translations.

 The nonramified part of this problem is the study of the space
$$
\cR^{nr}=\{f\colon T(F)\backslash{\bar E}/T(\cO)\to\Qlb\mid f(yz)=\zeta(y)f(z)\;\mbox{for all}\; y\in k^*/(k^*)^n\}
$$
as a representation of the corresponding Hecke algebra. For each $x\in X$ we get the Hecke algebra
\begin{multline*}
\cH_x=\{h\colon T(\cO)\backslash {\bar E_x}/T(\cO)\to\Qlb\mid h(yz)=\zeta(y)h(z)\;\mbox{for all}\; y\in k^*/(k^*)^n,\\
h \; \mbox{is of compact support}\; \}
\end{multline*}
with respect to convolution. Namely, if $h_i\in \cH_x$ then $h_1\ast h_2\in\cH_x$ is given as follows. For $u, z\in {\bar E_x}$ the expression $h_1(u)h_2(zu^{-1})$ depends only on the image of $u$ in $T(F_x)$, and we may define in this sense
$$
(h_1\ast h_2)(z)=\int_{u\in T(F_x)} h_1(u)h_2(zu^{-1})du,
$$
where $du$ is the Haar measure on $T(F_x)$ such that the volume of $T(\cO_x)$ is one. This algebra acts on $\cR^{nr}$ so that $h\in \cH_x$ acts on $f\in \cR^{nr}$ as
\begin{equation}
\label{def_action_h_on_f}
(f\ast h)(z)=\int_{u\in T(F_x)} f(zu^{-1})h(u)du.
\end{equation}
Again, it is understood that actually $u\in {\bar E_x}$, and the expression $f(zu^{-1})h(u)$ depends only on the image of $u$ in $T(F_x)$. With this definition, $\cR^{nr}$ is a left $\cH_x$-module: $(f\ast h_1)\ast h_2=f\ast (h_2\ast h_1)$. 

 For each $x\in X$ pick a uniformizer $t_x\in F_x$. Write $t_x^{\Lambda}$ for the image of the map $\Lambda\to T(F_x)$, $\lambda\mapsto t_x^{\lambda}$. Denote by $\Div(X,\Lambda)$ the group of $\Lambda$-valued divisors on $X$ viewed as a subgroup of $T(\AA)$. Namely, to a divisor $\sum_x \lambda_x x$ we associate $\prod_x t_x^{\lambda_x}\in T(\AA)$.

For the rest of Section~\ref{sec_classical} we make a stronger assumption:
\begin{itemize}
\item[(A)] the field $k$ satisfies $-1\in (k^*)^n$. 
\end{itemize} 
Since we are interested only in geometrizing the classical picture, this case is sufficient. This assumption implies that the restriction of (\ref{ext_BnEx_T(AA)_by_mu_n}) to $\Lambda$ is abelian. Pick a section of (\ref{ext_BnEx_T(AA)_by_mu_n}) over $\Lambda$. For $\lambda\in\Lambda$ we denote the corresponding element of $\bar E_x$ over $t_x^{\lambda}$ by $\bar t_x^{\lambda}$. 
 
  Note that for each $\lambda\in\Lambda$ there could be at most one function $h_{\lambda}\in \cH_x$ with support being the preimage of $T(\cO_x)t_x^{\lambda}T(\cO_x)$ in $\bar E_x$, and satisfying $h_{\lambda}(\bar t_x^{\lambda})=1$.
 
\begin{Lm} Let $\lambda\in\Lambda$. The following are equivalent.\\
i) $\lambda\in\Lambda^{\sharp}$,\\
ii) There is $h_{\lambda}\in \cH_x$ with support being the preimage of $T(\cO_x)t_x^{\lambda}T(\cO_x)$ in $\bar E_x$ and satisfying $h_{\lambda}(\bar t_x^{\lambda})=1$.
\end{Lm}
\begin{Prf}
For $v\in T(\cO_x)$ we have $v \bar t_x^{\lambda}=(v, t_x^{\lambda})_c \bar t_x^{\lambda} v$ in $\bar E_x$. So, $h_{\lambda}$ exists iff $(v, t_x^{\lambda})_c\in (k^*)^n$ for all $v\in T(\cO_x)$. For $v=\mu\otimes u$ with $\mu\in\Lambda, u\in \cO^*_x$ we get $(\mu\otimes u, t_x^{\lambda})_c=\bar u^{-\kappa(\mu, \lambda)}$, where $\bar u\in k^*$ is the image of $u$. Finally, $\bar u^{\kappa(\mu, \lambda)}\in (k^*)^n$ for all $\bar u\in k^*$, $\mu\in\Lambda$ if and only if $\lambda\in\Lambda^{\sharp}$.
\end{Prf} 

\medskip

 Let $\check{T}^{\sharp}$ be the Langlands dual torus over $\Qlb$ to $T^{\sharp}$. Write $\Rep(\check{T}^{\sharp})$ for the category of finite-dimensional representations of $\check{T}^{\sharp}$ over $\Qlb$.
 
\begin{Pp} 
\label{Pp_Satake}
There is an isomorphism of rings $\cH_x\iso \Qlb[\Lambda^{\sharp}]\iso K_0(\Rep(\check{T}^{\sharp}))$, where $K_0$ is the Grothendieck group of the corresponding category.
\end{Pp}
\begin{Prf}
Let $Z_x^{\lambda}$ be the preimage of $T(\cO_x)t_x^{\lambda}T(\cO_x)$ in $\bar E_x$. Clearly, if $\lambda,\mu\in \Lambda^{\sharp}$ then
$h_{\lambda}\ast h_{\mu}$ has support included in $Z_x^{\lambda+\mu}$. It is easy to check that actually $h_{\lambda}\ast h_{\mu}=h_{\lambda+\mu}$. So, we get an isomorphism $\cH_x\iso \Qlb[\Lambda^{\sharp}]$.
\end{Prf} 
 
\medskip

 The isogeny $T^{\sharp}\to T$ is a surjective morphism of sheaves in \'etale topology on $\Sch/k$. Let $K$ denote its kernel. This is a finite group that fits into an exact sequence
$$
1\to K(F)\to T^{\sharp}(F)\to T(F)\to \H^1(\Spec F, K)\to 1\ .
$$   

 The center $Z(\bar E)$ acts on $T(F)\backslash{\bar E}/T(\cO)$ by multiplication. Let $\Bun_{T^{\sharp}}$ denote the stack of $T^{\sharp}$-torsors over $X$. We have the push forward map $i_X\colon \Bun_{T^{\sharp}}\to \Bun_T$. Define $\bar K$ as the cokernel of $\Bun_{T^{\sharp}}(k)\to \Bun_T(k)$, this is a finite subgroup $\bar K\subset \H^2(X, K)$.
 
\begin{Rem} 
\label{Rem_Z(barE)-orbits}
The set of $Z(\bar E)$-orbits on $T(F)\backslash{\bar E}/T(\cO)$ is the group $T(\AA)/(T(F)Z^{\dag}T(\cO))$, it identifies with $\Coker(T^{\sharp}(\AA)\to T(F)\backslash T(\AA)/T(\cO))$, that is, with $\bar K$. The group $\bar K$ is usually nontrivial. 
\end{Rem} 
 
\subsubsection{} 
  
\begin{Lm} (1) The group $t_x^{\Lambda}Z(\bar E_x)$ is a maximal abelian subgroup of $\bar E_x$. 

\smallskip\noindent
(2) The group $\Div(X,\Lambda)Z(\bar E)$ is a maximal abelian subgroup in $\bar E$.
\end{Lm}
\begin{Prf}
1) By (A), this group is abelian. Pick bases $(\epsilon_i), (\eta_j)$ of $\Lambda$ and define $d_i, e_i\in\ZZ$ as in the proof of Lemma~\ref{Lm_global_max_abelain_subgroup}. Then 
$$
Z^{\dag}_x=\{\prod_i (\epsilon_i\otimes g_i)\in T(F_x)\mid g_i\in (F_x^*)^{e_i}\; \mbox{for all}\; i\} \ .
$$ 
Let now $v=\prod_i \epsilon_i\otimes v_i\in T(\cO_x)$ with $v_i\in \cO^*_x$. Let $\bar v_i$ be the image of $v_i$ under $\cO^*_x\to k^*$.  Assume that $(t_x^{\lambda}, v)_c=1$ for all $\lambda\in\Lambda$. So, for each $1\le j\le r$ we get 
$$
(\prod_i \epsilon_i\otimes v_i, t_x^{\eta_j})_c=(v_j, t_x)_{st}^{d_j}\in (k^*)^n \ .
$$
So, $\bar v_j^{d_j}\in (k^*)^n$. By Remark~\ref{Rem_about_finite_field}, we get $\bar v_j\in (k^*)^{e_j}$ for all $j$. So, $v\in Z^{\dag}_x$.\\
2) This follows from 1).
\end{Prf}  
  
\subsubsection{} 
\label{Section_325}
Let $\chi\colon Z(\bar E)\to\Qlb^*$ be any continuous character trivial on $Z(\bar E)\cap T(\cO)$ and extending $\zeta\colon k^*/(k^*)^n\to\Qlb^*$. Pick any extension of $\chi$ to a character 
$$
\bar\chi\colon \Div(X,\Lambda)Z(\bar E)\to\Qlb^*, 
$$
such an extension always exists, since $Z(\bar E)$ is an open subgroup of $\Div(X,\Lambda)Z(\bar E)$. By restriction, it yields a character $\bar\chi_x\colon t_x^{\Lambda}Z(\bar E_x)\to \Qlb^*$ for each $x\in X$. Let 
$$
\pi_x=\{ h\colon  {\bar E_x}\to\Qlb\mid h(yz)=\bar\chi_x(y) h(z)\;\mbox{for all}\; 
y\in t_x^{\Lambda}Z(\bar E_x)\}
$$
be the induced representation of $\bar E_x$. It is irreducible (see \cite{W}, Theorem~3.1). The space $\pi_x$ is a `twisted version' of the space of functions on 
$$
\prod_{i=1}^rk^*/(k^*)^{e_i},
$$
where $e_i$ and $r$ are defined as in the proof of Lemma~\ref{Lm_global_max_abelain_subgroup}, so $\dim\pi_x=\prod_{i=1}^r e_i$. Since $t_x^{\Lambda}T(\cO_x)=T(F_x)$, we get $\dim(\pi_x^{T(\cO_x)})=1$. So, we form the restricted tensor product
$$
\pi={\mathop{\bigotimes}\limits_{x\in X}}'\pi_x
$$
with respect to the unique spherical vector $h_x\in\pi_x^{T(\cO_x)}$ satisfying $h_x(1)=1$. This is an irreducible representation   
of $\bar E$, and $\dim(\pi^{T(\cO)})=1$. Consider also
\begin{multline*}
\bar\pi=\{h\colon {\bar E}\to \Qlb\mid h(yz)=\bar\chi(y)h(z)\; \mbox{for all}\; y\in \Div(X,\Lambda)Z(\bar E),\\
\mbox{there is an open subgroup} \; \cK\subset T(\cO), h(zz_1)=h(z), z_1\in \cK\}\ .
\end{multline*}
This is a smooth representation of $\bar E$. We identify $\pi\iso\bar\pi$ via the map sending $\bigotimes_x h_x$ to $h\in\bar\pi$, where $h(z)=\prod_{x\in X} h_x(z_x)$. 

  Assume in addition that $\chi\colon Z(\bar E)\to\Qlb^*$ is trivial on $Z(\bar E)\cap T(F)$. In this case there is a unique character 
$$
\und{\chi}\colon T(F)Z(\bar E)\to\Qlb^*
$$ 
extending $\chi\colon Z(\bar E)\to\Qlb^*$ and trivial on $T(F)$. Consider the $\bar E$-submodule 
$$
\cR_{\chi}=\{f\in\cR\mid  f(yz)=\chi(y)f(z)\;\mbox{for all}\; y\in Z(\bar E)\}
$$
of $\cR$. The results of (\cite{KP}, Section~0.3) apply in this situation. Combining them with Lemma~\ref{Lm_global_max_abelain_subgroup} we obtain that $\cR_{\chi}$ is irreducible and $\bar\pi\;\iso\;\cR_{\chi}$ as $\bar E$-modules.


 Let us write down an explicit isomorphism $S\colon \bar\pi\iso \cR_{\chi}$ of $\bar E$-modules. According to (\cite{KP}, Section~0.3), changing if necessary the extension $\bar\chi$ of $\chi$, we may and do assume the following:
\begin{itemize} 
\item[(C):] $\bar\chi$ and $\und{\chi}$ coincide over
$$
(\Div(X,\Lambda)Z(\bar E))\cap (T(F)Z(\bar E))\ .
$$ 
\end{itemize}
Under this additional assumption the isomorphism $S\colon \bar\pi\iso \cR_{\chi}$ sends $h$ to $Sh$, where 
\begin{equation}
\label{def_operator_S}
(Sh)(z)=\sum_{y\in T(F)/(\Div(X,\Lambda)Z(\bar E))\cap T(F)} h(yz), 
\end{equation}
here $z\in {\bar E}$. According to (C), if $u\in (\Div(X,\Lambda)Z(\bar E))\cap T(F)$ then $\bar\chi(u)=1$. So, for $y\in T(F)$ the expression $h(yz)$ depends
only on the image of $y$ in 
\begin{equation}
\label{set_summation_for_S}
T(F)/(\Div(X,\Lambda)Z(\bar E))\cap T(F).
\end{equation}
By the results of (\cite{KP}, Section~0.3), the map $S$ is an isomorphism of $\bar E$-representations. 

 The natural map $T^{\sharp}(F)\to T(F)$ factors through $Z(\bar E)\cap T(F)=Z^{\dag}\cap T(F)$. From Grunwald-Wang Theorem (\cite{AT}, Theorem~1, Chapter~IX) we see that the natural map $T^{\sharp}(F)\to T(F)\cap Z^{\dag}$ is surjective. Thus,
$$
T(F)/(Z(\bar E)\cap T(F))\iso \H^1(\Spec F, K).
$$
This is a finite group. The group $T(F)/(\Div(X,\Lambda)Z(\bar E))\cap T(F)$ is a quotient of the finite group $\H^1(\Spec F, K)$, so the sum in (\ref{def_operator_S}) is finite.

 Let $\phi\in \bar\pi^{T(\cO)}$ be the unique function satisfying $\phi(1)=1$. In Section~\ref{section_Geomerization} we will geometrize the \select{theta-function} 
$$
S\phi\colon  T(F)\backslash{\bar E}/T(\cO)\to\Qlb \ .
$$
We will construct a local system on the corresponding $\mu_n$-gerbe over $\Bun_T$, whose trace of Frobenius is the function $S\phi$ (namely, the complex $\cK$ in Proposition~\ref{Pp_H_eigensheaves}). 

\begin{Rem} The group $\H^1(\Spec F, K)$ is equipped with a skew-symmetric pairing with values in $\mu_n(k)$, it comes from the cup-product on $\H^1(\Spec F, K)$ and the bilinear form $\kappa$. We think that 
$(\Div(X,\Lambda)Z(\bar E))\cap T(F)/(Z(\bar E)\cap T(F))$ is a maximal isotropic subgroup in $\H^1(\Spec F, K)$ with respect to this pairing, but we have not checked this.
\end{Rem}

\subsubsection{} For each $x\in X$ let $\chi_x\colon \Lambda^{\sharp}\to\Qlb^*$ be the character sending $\mu$ to $\chi(\bar t_x^{\mu})$. The character $\chi_x^{-1}$ extends uniquely to an algebra homomorphism $\chi_x^{-1}\colon \Qlb[\Lambda^{\sharp}]\to\Qlb$. The 1-dimensional space $\pi_x^{T(\cO_x)}$ is naturally a module over $\cH_x$, the action being defined as in (\ref{def_action_h_on_f}). Then $\cH_x$ acts on $\pi_x^{T(\cO_x)}$ via the character 
\begin{equation}
\label{def_character_of_cHx}
\cH_x\;\iso\; \Qlb[\Lambda^{\sharp}]\;\toup{\chi_x^{-1}}\;\Qlb,
\end{equation}
the first isomorphism being that of Proposition~\ref{Pp_Satake}. 

 Set $\cR^{nr}_{\chi}=\cR^{nr}\cap \cR_{\chi}$, recall that this space is 1-dimensional. On the other hand, the $Z(\bar E)$-orbits on $T(F)\backslash \bar E/T(\cO)$ are identified with the group $\bar K$ by Remark~\ref{Rem_Z(barE)-orbits}. So, only one $Z(\bar E)$-orbit on $T(F)\backslash{\bar E}/T(\cO)$ supports a nonzero function from $\cR^{nr}_{\chi}$, one checks that this is the orbit through 1. 
 
\subsubsection{} Let $\deg\colon \Div(X,\Lambda)\to\Lambda$ be the degree map sending $\prod_x t_x^{\mu_x}$ to $\sum_x \mu_x$. For $\mu\in\Lambda$ write $\Div(X,\Lambda)^{\mu}$ for the subgroup of divisors of degree $\mu$.  One has $T(\AA)=\Div(X,\Lambda) T(\cO)$. For $\mu\in\Lambda$ write $T(\AA)^{\mu}=\Div(X,\Lambda)^{\mu} T(\cO)$. Note that $T(F)\subset T(\AA)^0$. 

 For $\mu\in\Lambda$ write $\bar E^{\mu}$ for the preimage of $T(\AA)^{\mu}$ under the projection $\bar E\to T(\AA)$. Let $Z(\bar E)^0=Z(\bar E)\cap \bar E^0$. The group $Z(\bar E)^0$ acts naturally on $T(F)\backslash {\bar E^{\mu}}/T(\cO)$.
 
  For $\mu\in\Lambda$ let
$$
^{\mu}\cR^{nr}_{\chi}=\{f\colon T(F)\backslash {\bar E^{\mu}}/T(\cO)\to\Qlb\mid f(yz)=\chi(y)f(z)\;\mbox{for}\; y\in Z(\bar E)^0\} \ .
$$
The property $\dim\cR_{\chi}^{nr}=1$ implies the following. If $\mu\notin\Lambda^{\sharp}$ then $^{\mu}\cR^{nr}_{\chi}=0$. If $\mu\in\Lambda^{\sharp}$ then $\dim(^{\mu}\cR^{nr}_{\chi})=1$, and any nonzero function from $^{\mu}\cR^{nr}_{\chi}$ is supported at $T(F)Z(\bar E)^{\mu}T(\cO)$. Set also 
$$
^{\mu}\cR^{nr}=\{f: T(F)\backslash {\bar E^{\mu}}/T(\cO)\to\Qlb\mid f(yz)=\zeta(y)f(z)\;\mbox{for}\; y\in k^*/(k^*)^n\}.
$$
The twisted Langlands correspondence on the classical level becomes essentially the following.

\begin{Pp} Let $\mu\in\Lambda$. \\
(1) If $\mu\notin\Lambda^{\sharp}$ then $^{\mu}\cR^{nr}=0$.\\
(2) If $\mu\in\Lambda^{\sharp}$ then there is a finite direct sum decomposition
$$
^{\mu}\cR^{nr}=\bigoplus_{\chi} \;\; {^{\mu}\cR^{nr}_{\chi}},
$$
the sum over all characters $\chi\colon Z(\bar E)^0\to \Qlb^*$ trivial on $Z(\bar E)^0\cap T(\cO)$ and on $Z(\bar E)^0\cap T(F)$ and extending $\zeta\colon k^*/(k^*)^n\to\Qlb^*$. \QED
\end{Pp}
  
\section{Preliminaries to geometrization}
\label{section_Preliminaries}

\subsection{Gerbs via central extensions} 

\subsubsection{} From now on we assume $k$ algebraically closed of characteristic $p\ge 0$. Let $Z$ be a $k$-scheme, $n\ge 1$ invertible in $k$. Let $G$ be a finite group acting on $Z$. A lifting of this action on the trivial gerbe $Z\times B(\mu_n)$ is described as follows. 

 For $g\in G$ and a $\mu_n$-torsor $\cF$ on $Z$ we define a morphism $(\cF, g)\colon Z\times B(\mu_n)\to Z\times B(\mu_n)$. For a $S$-point $(f, \cT)$, where $f\colon S\to Z$ and $\cT$ is a $\mu_n$-torsor on $S$, the map $(\cF, g)$
sends it to the $S$-point $(gf, \cT\otimes f^*\cF)$. The composition $(\cF_2, g_2) (\cF_1, g_1)\iso (\cF_1\otimes g_1^*\cF_2, g_2g_1)$ canonically. 

 The action of $G$ on $Z\times B(\mu_n)$ is given by the data: for each $g\in G$ a $\mu_n$-torsor $\cF_g$ on $Z$. We assume that $\cF_1=\cF^0$ is the trivial  
$\mu_n$-torsor. For each pair $g,h\in G$ we are given an isomorphism $\tau_{g,h}\colon \cF_h\otimes h^*\cF_g\iso \cF_{gh}$ of $\mu_n$-torsors on $Z$. It is required that for any $g,h,x\in G$ the diagram commutes
$$
\begin{array}{ccc} 
\cF_x\otimes x^*\cF_h\otimes x^*h^*\cF_g & \toup{\tau_{h,x}} & \cF_{hx}\otimes x^*h^*\cF_g\\
\downarrow\lefteqn{\scriptstyle x^*\tau_{g,h}} && \downarrow\lefteqn{\scriptstyle \tau_{g, hx}}\\
\cF_x\otimes x^*\cF_{gh} & \toup{\tau_{gh, x}} & \cF_{ghx} \ .
\end{array}
$$
Besides, for $g\in G$ each of the maps $\tau_{g,1}\colon \cF_1\otimes \cF_g\to\cF_g$ and $\tau_{1,g}\colon \cF_g\otimes g^*\cF_1\to\cF_g$ is the identity. 

 Write $\Tors(Z,\mu_n)$ for the groupoid of $\mu_n$-torsors on $Z$. Consider the groupoid $G_Z$ of pairs $(\cF, g)$, where $\cF\in \Tors(Z,\mu_n)$, $g\in G$. We define the multiplication functor $m\colon G_Z\times G_Z\to G_Z$ on objects by 
$$
m((\cF_2, g_2) (\cF_1, g_1))=(\cF_1\otimes g_1^*\cF_2, g_2g_1).
$$
It is defined naturally on morphisms. This extends naturally to a structure of a group stack on the groupoid $G_Z$. We obtain an exact sequence of group stacks
\begin{equation}
\label{ext_G_by_Tors(Zmu_n)}
1\to \Tors(Z,\mu_n)\to G_Z\to G\to 1.
\end{equation}
The group stack $G_Z$ acts on $Z\times B(\mu_n)$, namely $(\cF, g)$ acts by the morphism $(\cF, g)$. So, the datum of an action of $G$ on $Z\times B(\mu_n)$ is a section of (\ref{ext_G_by_Tors(Zmu_n)}) in the category of group stacks. That is, a morphism of group stacks $G\to G_Z$ whose projection to $G$ is the identity. Such morphism always exists.

  Let $\und{G}_Z$ and $Tors(Z, \mu_n)$ be the coarse moduli space of $G_Z$ and $\Tors(Z,\mu_n)$ respectively. These are abstract groups. We get an exact sequence
$$
1\to Tors(Z,\mu_n)\to \und{G}_Z\to G\to 1,
$$
which is a semi-direct product with respect to the action of $G$ on $Tors(Z,\mu)$ such that $g\in G$ sends $\cF$ to $(g^{-1})^*\cF$. 

 If $\cT\in \Tors(Z,\mu_n)$ and a section of (\ref{ext_G_by_Tors(Zmu_n)}) is given by a collection $(\cF_g, \tau_{g,h})$, $g,h\in G$ as above, we may congugate this section by the element $(\cT,1)\in G_Z$. This produces the collection $(\cF'_g, \tau'_{g,h})$, where $\cF'_g=\cT^{-1}\otimes g^*\cT\otimes \cF_g$ and $\tau'=\tau$.

\begin{Rem}
\label{Rem_for_Section_B1}
(1) Now let $1\to \mu_n(k)\to\tilde G\to G\to 1$ be a central extension in the category of groups. Take $\cF_g$ be the constant $\mu_n$-torsor consisting of all $\tilde g\in  \tilde G$ over $g$. The group structure of $\tilde G$ yields an isomorphism $\tau_{g,h}\colon \cF_h\otimes\cF_g\iso \cF_{gh}$. The above conditions on $\tau$ are verified, so we get an action of $G$ on $Z\times B(\mu_n)$ extending the action of $G$ on $Z$. 
 
  If we let $\tilde G$ act on $Z$ via the homomorphism $\tilde G\to G$ with the previous action of $G$ then the stack quotient $Z/\tilde G$ is a $\mu_n$-gerbe over the stack quotient $Z/G$. 
  
\smallskip\noindent  
(2) Conversely, let $\cY\to Z/G$ be a $\mu_n$-gerbe equipped with a trivialization of the gerbe ${\cY\times_{Z/G} Z}\to Z$. This yields a section of (\ref{ext_G_by_Tors(Zmu_n)}) given by a collection $(\cF_g, \tau_{g,h})$, $g,h\in G$ as above. Assume also that all the $\mu_n$-torsors $\cF_g$ are trivial. So, we view $\cF_g$ as a $\mu_n$-torsor over a point. Then this $\mu_n$-torsor over $G$ is multiplicative in the sense of \cite{DB}, so this is a central extension $1\to \mu_n\to \tilde G\to G\to 1$. Moreover, the trivialization of $\cY\times_{Z/G} Z\to Z$ descends to an isomorphism $\cY\iso Z/\tilde G$ over $Z/G$. 
\end{Rem}
  
\subsubsection{} Assume that $\alpha\colon Z'\to Z$ is a Galois \'etale covering with Galois group $\Gamma$. Assume also that for each $g\in G$ there is $g'\in\Aut(Z')$ such that the diagram commutes
$$
\begin{array}{ccc}
Z' & \toup{\alpha} &Z\\
\downarrow\lefteqn{\scriptstyle g'} && \downarrow\lefteqn{\scriptstyle g}\\
Z' & \toup{\alpha} &Z .
\end{array}
$$
We get the group $G'$ of automorphisms $g'\in\Aut(Z')$ for which such $g\in G$ exists, it fits into an exact sequence $1\to \Gamma\to G'\toup{\beta} G\to 1$. 

 Assume a section of (\ref{ext_G_by_Tors(Zmu_n)}) is given by a collection $(\cF_g, \tau_{g,h})$, $g,h\in G$ as above, assume also that all the $\mu_n$-torsors $\alpha^*\cF_g$ are trivial over $Z'$. Consider the exact sequence
$$
1\to\Tors(Z',\mu_n)\to G'_Z\to G'\to 1
$$
defined for the the above action of $G'$ on $Z'$. We get a section of this exact sequence given by the collection $(g', \cF_{g'}, \tau')$, where $\cF_{g'}=\alpha^*\cF_g$ for $g=\beta(g')$, and $\tau'=\alpha^*\tau$. 

 By Remark~\ref{Rem_for_Section_B1} (2), we get a central extension 
\begin{equation}
\label{ext_for_B2} 
1\to \mu_n\to \tilde G'\to G'\to 1
\end{equation}
and a $\mu_n$-gerbe $Z'/\tilde G'\to Z/G$.   

\begin{Rem} One may ask if any $\mu_n$-gerbe over the stack quotient $Z/G$ comes from a central extension $1\to \mu_n\to ?\to G\to 1$. We will not answer this question in this paper, but we think the answer is `no'. If (\ref{ext_for_B2}) does not admit a section over $\Gamma$, whose image is a normal subgroup in $\tilde G'$ then the gerbe $Z'/\tilde G'$ would provide a counterexample.
\end{Rem}

\subsection{$\theta$-data and central extensions of $\Bun_T$}
\label{subsection_theta_data}

\subsubsection{} 
\label{section_421}
Let $\Lambda$ be a free abelian group of finite type, $\check{\Lambda}=\Hom(\Lambda,\ZZ)$. Set $T=\Lambda\otimes\Gm$ and $\check{T}=\check{\Lambda}\otimes\Gm$. Let $X$ be a smooth projective connected curve over $k$. 

We will use Picard groupoid $\cP^{\theta}(X,\Lambda)$ of $\theta$-data introduced by Beilinson-Drinfeld in (\cite{BD2}, Section~3.10.3). Recall that an object of $\cP^{\theta}(X,\Lambda)$ is a triple $\theta=(\kappa,\lambda, c)$, where $\kappa\colon \Lambda\otimes\Lambda\to\ZZ$ is a symmetric bilinear form, $\lambda$ is a rule that assigns to each $\gamma\in\Lambda$ a super line bundle $\lambda^{\gamma}$ on $X$, and $c$ is a rule that assigns to each pair $\gamma_1,\gamma_2\in\Lambda$ an isomorphism $c^{\gamma_1,\gamma_2}\colon  \lambda^{\gamma_1}\otimes\lambda^{\gamma_2}\iso \lambda^{\gamma_1+\gamma_2}\otimes\Omega^{\kappa(\gamma_1,\gamma_2)}$ on $X$. They are subject to the conditions explained in (\cite{BD2}, Section~3.10.3). In particular, recall that the parity of $\lambda^{\gamma}$ is $\kappa(\gamma,\gamma)\!\!\mod 2$.

 Fixing a symmetric bilinear form $\kappa\colon \Lambda\otimes\Lambda\to\ZZ$, one gets a subgroupoid $\cP^{\theta}(X,\Lambda)^{\kappa}\subset \cP^{\theta}(X,\Lambda)$. Recall that $\cP(X,\Lambda):=\cP^{\theta}(X,\Lambda)^0$ is a Picard subgroupoid, and each $\cP^{\theta}(X,\Lambda)^{\kappa}$ is a $\cP(X,\Lambda)$-torsor. By (\cite{BD2}, 3.10.3.1), there is a canonical equivalence of Picard groupoids
\begin{equation}
\label{equiv_cP^theta(X,Lambda)}
\cP(X,\Lambda)\iso \Tors(X, \check{T}),
\end{equation}
where $\Tors(X, \check{T})$ is the Picard groupoid of $\check{T}$-torsors on $X$.
Recall the groupoid $\cE^s(T)$ defined in Section~\ref{section_311}. The following construction is borrowed from (\cite{BD2}, Lemma~3.10.3.1).

\begin{Lm} 
\label{Lm_functor_from_Weissman_to_theta_data}
Pick a square root $\Omega^{\frac{1}{2}}$ of $\Omega$ on $X$. It gives rise to a functor $\cE^s(T)\to \cP^{\theta}(X,\Lambda)$.
\end{Lm}
\begin{Prf}
Let $(\kappa, \tilde\Lambda^s)\in \cE^s(T)$, so for each $\gamma\in\Lambda$ we are given a super line $\epsilon^{\gamma}$ and isomorphisms $c^{\gamma_1,\gamma_2}\colon \epsilon^{\gamma_1}\otimes \epsilon^{\gamma_2}\iso \epsilon^{\gamma_1+\gamma_2}$. For $\gamma\in\Lambda$ set $\lambda^{\gamma}=(\Omega^{\frac{1}{2}})^{\otimes -\kappa(\gamma,\gamma)}\otimes \epsilon^{\gamma}$. Let $'c^{\gamma_1,\gamma_2}\colon \lambda^{\gamma_1}\otimes\lambda^{\gamma_2}\iso \lambda^{\gamma_1+\gamma_2}\otimes\Omega^{\kappa(\gamma_1,\gamma_2)}$ be the evident product obtained from $c^{\gamma_1,\gamma_2}$. Then $(\kappa,\lambda, 'c)\in\cP^{\theta}(X,\Lambda)$. 
\end{Prf}

\medskip
 
 One has the sheaf $\Div(X,\Lambda)$ on the category $\Sch/k$ of $k$-schemes in flat topology introduced in (\cite{BD2}, 3.10.7). It classifies relative $\Lambda$-valued Cartier divisors on $X$. One has the Abel-Jacobi map $AJ\colon \Div(X,\Lambda)\to\Bun_T$ given by $D\otimes \gamma\mapsto \cO(D)^{\otimes\gamma}$ for $\gamma\in\Lambda$, $D\in\Div(X,\ZZ)$. This is a morphism of abelian group stacks. 
 
 In (\cite{BD2}, Section~3.10.7) the Picard groupoid $\cPic^f(\Div(X,\Lambda))$ of factorizable super line bundles on $\Div(X,\Lambda)$ is introduced. By (\cite{BD2}, Proposition~3.10.7.1), one has a natural equivalence of Picard groupoids
\begin{equation}
\label{equivalence_cPic^f(Div(X,Lambda)_and_cP^theta(X,Lambda)}
\cPic^f(\Div(X,\Lambda))\iso \cP^{\theta}(X,\Lambda).
\end{equation}
Write $\cPic(\Bun_T)$ and $\cPic(\Div(X,\Lambda))$ for the Picard groupoids of super line bundles on $\Bun_T$ and $\Div(X,\Lambda)$ respectively. By (\cite{BD2}, Proposition~4.9.1.2), the functor $AJ^*$ is an equivalence of Picard groupoids 
$$
\cPic(\Bun_T)\iso \cPic(\Div(X,\Lambda))\ .
$$

For $\mu\in\Lambda$ write $\Bun_T^{\mu}$ for the connected component of $\Bun_T$ classifying $\cF\in\Bun_T$ such that for any $\check{\lambda}\in\check{\Lambda}$ one has $\deg(\cL^{\check{\lambda}}_{\cF})=\<\mu, \check{\lambda}\>$. 

 Recall that each line bundle $\tau$ on $\Bun_T$ defines a map $\delta_{\tau}\colon \Lambda\to\check{\Lambda}$ such that for $\cF\in \Bun_T^{\mu}$ the group $T(k)\subset \Aut(\cF)$ acts on the fibre $\tau_{\cF}$ by $\delta_{\tau}(\mu)$. 

The forgetful functor $\cPic^f(\Div(X,\Lambda))\to \cPic(\Div(X,\Lambda))$ yields a composition
\begin{equation}
\label{functor_from_cP^theta(X,Lambda)_to_cPic(Bun_T)}
\cP^{\theta}(X,\Lambda)\iso \cPic^f(\Div(X,\Lambda))\to \cPic(\Bun_T).
\end{equation}

\subsubsection{}
\label{section_422}
 For a $\check{T}$-torsor $\cT$ on $X$ let $L_{\cT}$ denote the factorizable line bundle on $\Div(X,\Lambda)$ associated to $\cT$ via (\ref{equiv_cP^theta(X,Lambda)}) and 
(\ref{equivalence_cPic^f(Div(X,Lambda)_and_cP^theta(X,Lambda)}), it is of parity zero. For a $k$-point $D=\sum_x \lambda_x x$ of $\Div(X,\Lambda)$ with $\lambda_x\in\Lambda$ the fibre of $L_{\cT}$ at $D$ is $(L_{\cT})_D=\bigotimes_{x\in X} (\cL^{\lambda_x}_{\cT})_x$. For more general points the construction of $(L_{\cT})_D$ is based on the norm map (see the proof of \cite{BD2}, Proposition~3.10.7.1). 

 As $\cT$ varies in $\Bun_{\check{T}}$, these line bundles form a line bundle $L$ on $\Bun_{\check{T}}\times\Div(X,\Lambda)$. As in (\cite{BD2}, Proposition~4.9.1.2), one checks that 
there is a line bundle $\cL^{univ}$ on $\Bun_{\check{T}}\times\Bun_T$ equipped with an isomorphism $(\id\times AJ)^*\cL^{univ}\iso L$, where 
\begin{equation}
\label{map_q_appendix_C}
\id\times AJ\colon \Bun_{\check{T}}\times\Div(X,\Lambda)\to \Bun_{\check{T}}\times\Bun_T.
\end{equation}
The line bundle $\cL^{univ}$ is defined up to a unique isomorphism. 
  
   Let $\cL^{univ}_{\cT}$ denote the restriction of $\cL^{univ}$ to $\Bun_T$ given by fixing a $k$-point $\cT$ of $\Bun_{\check{T}}$. By (\cite{BD2}, Lemma~4.9.2), the map $\delta_{\cL^{univ}_{\cT}}$ is constant, its image equals $\deg(\cT)\in \check{\Lambda}$. 
   
   The line bundle $\cL^{univ}$ defines a biextension of $\Bun_{\check{T}}\times\Bun_T$ in the sense of (\cite{P}, Section~10.3). So, $\cL^{univ}$ can be seen as a commutative central extension of $\Bun_{\check{T}}\times\Bun_T$ by $\Gm\times \Bun_{\check{T}}$ in the category of commutative group stacks over $\Bun_{\check{T}}$, and also as a commutative central extension of $\Bun_{\check{T}}\times\Bun_T$ by $\Gm\times \Bun_{T}$ in the category of commutative group stacks over $\Bun_{T}$.
 
  For $n\ge 1$ write $\Bun_n$ for the stack of rank $n$ vector bundles on $X$.
For example, if $T=\Gm$ then $\check{T}=\Gm$, and the line bundle $\cL^{univ}$ over $\Bun_1\times\Bun_1$ is canonically isomorphic to the line bundle, whose fibre at $(\cA_1, \cA_2)$ is 
$$
\frac{\det\RG(X, \cA_1\otimes\cA_2)\otimes\det\RG(X,\cO)}{\det\RG(X,\cA_1)\otimes\det\RG(X,\cA_2)} \ .
$$
\begin{Rem} 
\label{Rem_some_symmetry}
(1) The line bundle $\cL^{univ}$ is symmetric in the following sense. We could start with a $T$-torsor $\cF$ on $X$ and consider the corresponding factorizable line bundle $L_{\cF}$ on $\Div(X,\check{\Lambda})$. As $\cF$ varies, they form a line bundle $L$ on $\Div(X, \check{\Lambda})\times\Bun_T$. Then $L$ is canonically isomorphic to $(AJ\times\id)^*\cL^{univ}$, where
$$
AJ\times\id\colon \Div(X, \check{\Lambda})\times\Bun_T\to \Bun_{\check{T}}\times\Bun_T .
$$
(2) More generally, if $\Lambda'$ is a free abelian group of finite type, let $\eta\colon \Lambda\times\Lambda'\to\ZZ$ be a bilinear form. Let $T'=\Lambda'\otimes\Gm$. The corresponding linear maps $\eta\colon \Lambda\to \check{\Lambda}'$ and $\eta\colon \Lambda'\to \check{\Lambda}$ yield maps $T\to \check{T}'$ and $T'\to \check{T}$ respectively, hence a diagram
$$
\Bun_{\check{T}'}\times\Bun_{T'}\getsup{\eta_X\times\id} \Bun_T\times\Bun_{T'} \toup{\id\times\eta_X} \Bun_T\times \Bun_{\check{T}} .
$$
For this diagram the biextensions $(\eta_X\times\id)^*\cL^{univ}$ and $(\id\times\eta_X)^*\cL^{univ}$ of $\Bun_T\times\Bun_{T'}$ are canonically isomorphic, we denote this biextension by $^{\eta}\cL^{univ}$. 
\end{Rem}

\subsubsection{} 
\label{section_423}
For $\check{\lambda}\in\check{\Lambda}$ let $L^{\check{\lambda}}$ denote the line bundle on $\Bun_T$ with fibre 
$$
\det\RG(X, \cL^{\check{\lambda}}_{\cF})\otimes \det\RG(X,\cO)^{-1}
$$ 
at $\cF\in\Bun_T$. We view it as $\ZZ/2\ZZ$-graded of parity $\<\mu,\check{\lambda}\>$ over $\Bun_T^{\mu}$. Let $\kappa\colon\Lambda\otimes\Lambda\to\ZZ$ be a symmetric bilinear form. Given a presentation denoted $\beta$
$$
\kappa=\sum_i b_i(\check{\lambda}_i\otimes\check{\lambda}_i)
$$
for some $b_i\in\ZZ$, $\check{\lambda}_i\in\check{\Lambda}$ 
we associate to it a line bundle $L_{\beta}=\bigotimes_i (L^{\check{\lambda}_i})^{\otimes b_i}$ on $\Bun_T$. This is the image of some element of $\cP^{\theta}(X,\Lambda)^{\kappa}$ under (\ref{functor_from_cP^theta(X,Lambda)_to_cPic(Bun_T)}), compare with Lemma~\ref{Lm_1_realizations}. 

 Let $\theta=(\lambda, \kappa,c)\in \cP^{\theta}(X,\Lambda)$, write also $\lambda$ for its image under (\ref{functor_from_cP^theta(X,Lambda)_to_cPic(Bun_T)}). 

\begin{Pp} 
\label{Pp_factoriztion_for_twists_by_divisors}
For $x\in X$, $\mu\in \Lambda$, $\cF\in\Bun_T$ there is a natural $\ZZ/2\ZZ$-graded isomorphism
\begin{equation}
\label{iso_we_do_need}
\lambda_{\cF(\mu x)}\iso \lambda_{\cF} \otimes (\cL^{\kappa(\mu)}_{\cF})_x\otimes\lambda_{\cO(\mu x)}
\end{equation}
functorial in $\theta\in \cP^{\theta}(X,\Lambda)$. Here $(\cL^{\kappa(\mu)}_{\cF})_x$ is of parity zero.
\end{Pp}
\begin{Prf}
1) First, consider the case $\kappa=0$. Recall that the image of an object of $\cP(X,\Lambda)$ under (\ref{functor_from_cP^theta(X,Lambda)_to_cPic(Bun_T)}) is a central extension of $\Bun_T$ by $\Gm$ (cf. \cite{BD2}, Lemma~4.9.2). In this case the line bundle $\lambda$ on $\Bun_T$ is multiplicative. So, for any $\cF,\cF'\in\Bun_T$ one has $\lambda_{\cF\otimes\cF'}\iso \lambda_{\cF}\otimes\lambda_{\cF'}$ and we are done.

\smallskip
\noindent
2) Assume first that $\kappa=\check{\lambda}\otimes\check{\lambda}$ and $\lambda$ is of the form $L^{\check{\lambda}}$ as above, it comes from some particular object of $\cP^{\theta}(X,\Lambda)^{\kappa}$. In this case $\kappa(\mu)=\<\mu, \check{\lambda}\>\check{\lambda}$. Assume also $\<\mu, \check{\lambda}\>\ge 0$, the opposite case being similar. Then 
$$
L^{\check{\lambda}}_{\cF(\mu x)}\,\iso\, L^{\check{\lambda}}_{\cF}\otimes \det\RG(X, \cL^{\check{\lambda}}_{\cF}(\<\mu, \check{\lambda}\>x)/\cL^{\check{\lambda}}_{\cF})
$$
and $L^{\check{\lambda}}_{\cO(\mu x)}\,\iso\, \det\RG(X, \cO(\<\mu, \check{\lambda}\>x)/\cO)$ canonically. To conclude, note that for a line bundle $\cA$ on $X$, $m\ge 0$, one has $\det\RG(X, \cA(mx)/\cA)\,\iso\, \cA_x^m\otimes \det\RG(X, \cO(mx)/\cO)$ canonically.

\smallskip
\noindent
3) For general $\kappa$ pick a presentation $\kappa=\sum_i b_i(\check{\lambda}\otimes \check{\lambda})$, it gives rise to the line bundle $L_{\beta}$ on $\Bun_T$ coming from some particular object of $\cP^{\theta}(X,\Lambda)^{\kappa}$. The desired isomorphism in this case is the product of isomorphisms obtained in 2).

 Moreover, the isomorphism (\ref{iso_we_do_need}) is equivariant with respect to the action of $\check{T}$, here $\check{T}$ is the group of automorphisms of any object of $\cP^{\theta}(X,\Lambda)$. Indeed, $z\in\check{T}$ acts on $\lambda\in\Pic(\Bun_T)$ so that it acts on the fibre $\lambda_{\cF}$ as $(\deg\cF)(z)$. Combining with 1), we get the desired isomorphism in general. 
\end{Prf}

\medskip

 Write $\bar\kappa\colon \Bun_T\to\Bun_{\check{T}}$ for the map sending $\cF$ to the $\check{T}$-torsor $\bar\kappa(\cF)$ such that for $\lambda\in\Lambda$ one has $\cL^{\lambda}_{\bar\kappa(\cF)}=\cL^{\kappa(\lambda)}_{\cF}$. 
 
\begin{Pp} 
\label{Pp7_factorization_of_lambda_on_BunT}
For $\cF,\cT\in\Bun_T$ there is a natural $\ZZ/2\ZZ$-graded isomorphism
$$
\lambda_{\cF\otimes\cT}\iso \lambda_{\cF}\otimes\lambda_{\cT}\otimes {^{\kappa}\!\cL^{univ}_{\cF, \cT}} \ .
$$
Here $^{\kappa}\!\cL^{univ}$ is the line bundle purely of parity zero on $\Bun_{T}\times\Bun_T$ defined in Remark~\ref{Rem_some_symmetry}.
\end{Pp} 
\begin{Prf} Consider the map $\id\times AJ\colon \Bun_T\times\Div(X,\Lambda)\to\Bun_T\times\Bun_T$. By (\cite{BD}, Proposition~4.9.1.2), is suffices to establish the desired isomorphism after applying $(\id\times AJ)^*$ to both sides. So, we may assume $\cT=\cO(D)$, where $D\in \Div(X,\Lambda)$. If $D=\sum_i \mu_i x_i$ with $x_i$ pairwise distinct then Proposition~\ref{Pp_factoriztion_for_twists_by_divisors} gives an isomorphism
$$
\lambda_{\cF(D)}\iso\lambda_{\cF}\otimes\lambda_{\cO(D)}\otimes (\otimes_i (\cL^{\kappa(\mu_i)}_{\cF})_{x_i}) \ .
$$
With the notations of Section~\ref{section_422}, we get $(\otimes_i (\cL^{\kappa(\mu_i)}_{\cF})_{x_i})\iso (L_{\bar\kappa(\cF)})_D$, the fibre of $L_{\bar\kappa(\cF)}\in\Pic(\Div(X,\Lambda))$ at $D$. The latter identifies canonically with $\cL^{univ}_{\bar\kappa(\cF), \cT}\iso {^{\kappa}\cL^{univ}_{\cF, \cT}}$. Our claim follows.
\end{Prf} 
 
\begin{Rem} For $\mu'\in \Lambda$ let $\cF=\cO(\mu'x)$. Then the isomorphism (\ref{iso_we_do_need}) becomes the isomorphism $c^{\mu, \mu'}$ in the definition of the theta-datum. 
\end{Rem}

\section{Geometrization}
\label{section_Geomerization}

\subsection{Local setting}

\subsubsection{} Let $\Lambda, \check{\Lambda}, T$ be as in Section~\ref{section_421}. Let $\cO=k[[t]]\subset F=k((t))$. View $T(F)$ as a commutative group ind-scheme over $k$. Recall that Contou-Carr\`ere defined in \cite{CC} a canonical skew-symmetric symbol $(\cdot, \cdot)_{st}: F^*\times F^*\to\Gm$. This is a morphism of ind-schemes, on the level of $k$-points equal to the tame symbol (see \cite{BBE}, Sections~3.1-3.3). 

 Let $\kappa\colon \Lambda\otimes\Lambda\to\ZZ$ be an even symmetric bilinear form. In (\cite{BD}, Definition~10.3.13) a notion of a Heisenberg $\kappa$-extension of $T(F)$ was introduced. This is a central extension 
\begin{equation}
\label{ext_T(F)_by_Gm} 
 1\to\Gm\to \cE\to T(F)\to 1
\end{equation} 
in the category of group ind-schemes, whose commutator satisfies
$$
(\lambda_1\otimes f_1, \lambda_2\otimes f_2)_c=(f_1,f_2)_{st}^{-\kappa(\lambda_1,\lambda_2)}
$$ 
for $f_i\in F^*$, $\lambda_i\in \Lambda$. The Heisenberg extensions of $T(F)$ form a Picard groupoid, whose structure is described in (\cite{BD}, Section~10.3.13). Let (\ref{ext_T(F)_by_Gm}) be a Heisenberg $(-\kappa)$-extension of $T(F)$. 

  Pick $n>1$ invertible in $k$, pick a primitive character $\zeta\colon \mu_n(k)\to\Qlb^*$. Let $\Gm$ act on $\cE$ so that $y\in\Gm$ acts on $e\in E$ as $y^ne$. The stack quotient $E=\cE/\Gm$ under this action fits into an exact sequence of group stacks
\begin{equation}
\label{ext_T(F)_by_B(mu_n)}
1\to B(\mu_n)\to E \to T(F)\to 1.
\end{equation}
This is a geometric analog of the extension (\ref{ext_BnEx_T(AA)_by_mu_n}). 

 Let $\Lambda^{\sharp}$ be given by (\ref{def_Lambda^sharp}), so $\Lambda^{\sharp}\subset\Lambda$ is of finite index. Set $T^{\sharp}=\Lambda^{\sharp}\otimes\Gm$. Let $i\colon T^{\sharp}\to T$ be the corresponding isogeny.
 
 Let $E^{\sharp}$ (resp., $\cE^{\sharp}$) be obtained from $E$ (resp., $\cE$) by the base change $T^{\sharp}(F)\to T(F)$. The group stack $E^{\sharp}$ admits a natural commutativity constraint, so it is a commutative group stack. Moreover, consider the maps $m\colon E^{\sharp}\times E\to E$ and $m'\colon E^{\sharp}\times E\to E$, where $m$ is the product, and $m'(xy)=m(yx)$. 
 
\begin{Lm} 
\label{Lm_local_geom_one}
There is a natural 2-isomorphism $m'\iso m$. In this sense $E^{\sharp}$ is contained in the `center' of $E$.
\end{Lm}
\begin{Prf}
One has $T(F)\iso \Lambda\times(\Lambda\otimes (\Gm\times \WW\times \hat\WW))$ canonically, where $\WW$ is the group scheme of big Witt vectors, and $\hat\WW$ is its completion (with the change of variables $t\mapsto t^{-1}$), see (\cite{BBE}, Section~3.2). So, $T^{\sharp}(F)=\Lambda^{\sharp}\times(\Lambda^{\sharp}\otimes (\Gm\times \WW\times \hat\WW))$. The composition 
$T^{\sharp}(F)\times T(F)\to T(F)\times T(F)\toup{(\cdot,\cdot)_c}\Gm$ factors as 
$
T^{\sharp}(F)\times T(F)\to \Gm\toup{z\mapsto z^n} \Gm
$.    
\end{Prf}  
 
\begin{Rem} It is natural to ask for a description of the Drinfeld center of $E$. Is it an algebraic stack? What is its relation with $E^{\sharp}$? We will not need to answer these questions in this paper.
\end{Rem}  
 
\smallskip

  Set $T(F)'=\Lambda\times(\Lambda^{\sharp}\otimes (\Gm\times \WW\times \hat\WW))$ and 
$$
T(F)''=\Lambda^{\sharp}\times(\Lambda\otimes (\Gm\times \WW))\hook{} T(F)_{red} \ .
$$ 
Let $E'$ (resp., $E''$) be the base change of $E$ by $T(F)'\hook{} T(F)$ (resp., by $T(F)''\hook{} T(F)$). 
  
\begin{Lm}
(1) Assume that either $n$ is odd or $\kappa(\Lambda\otimes\Lambda)\subset 2\ZZ$. Then $E'$ is naturally an abelian group stack. \\
(2) If $n$ is even then $E''$ is naturally an abelian group stack. 
\end{Lm}
\begin{Prf}
(1) Consider the extension $1\to \Gm\to \tilde\Lambda\to \Lambda\to 1$ obtained as the pull-back of (\ref{ext_T(F)_by_Gm}) by $\Lambda\to T(F)$, $\lambda\mapsto t^{\lambda}$. It commutator is $(\lambda_1,\lambda_2)_c=(-1)^{\kappa(\lambda_1,\lambda_2)}$. If $\kappa(\Lambda\otimes\Lambda)\subset 2\ZZ$ then $\tilde\Lambda$ is abelian, so one gets a commutativity constraint for $E'$ as in Lemma~\ref{Lm_local_geom_one}. 

 We may view $\tilde\Lambda$ as an extension of $\Lambda$ by $\mu_2$ with the same commutator. We have a morphism of group stacks $\bar\delta\colon \mu_2\to B(\mu_n)$ given by the exact sequence of groups 
\begin{equation}
\label{seq_mu_2_by_mu_n} 
 1\to \mu_n\to \mu_{2n}\toup{\delta}\mu_2\to 1,
\end{equation}  
here $\delta(z)=z^n$. Namely, $\delta$ is a $\mu_n$-torsor, which can be viewed as a morphism $\bar\delta\colon \mu_2\to B(\mu_n)$. The multiplicative structure on the $\mu_n$-torsor $\delta$ provides a structure of a morphism of group stacks on $\bar\delta$. Let $1\to B(\mu_n)\to \bar\Lambda\to \Lambda\to 1$ be the push-forward of $1\to \mu_2\to \tilde\Lambda\to \Lambda\to 1$ by $\bar\delta$. Then $\bar\Lambda$ is the restriction of $E$ to $\Lambda$. 

If $n$ is odd then the sequence (\ref{seq_mu_2_by_mu_n}) splits canonically, so the restriction of $E$ to $\Lambda$ is naturally an abelian group stack. Now one gets the desired commutativity constraint for $E'$ as above.

\smallskip\noindent
(2) Let $\tilde\Lambda^{\sharp}$ denote the restriction of $\tilde\Lambda$ to $\Lambda^{\sharp}$. If $n$ is even then $\tilde\Lambda^{\sharp}$ is an abelian group scheme. Now we can construct a commutativity constraint on $E''$ as in 1) using the fact that the symbol $(\cdot,\cdot)_{st}$ is trivial on $\Gm\times\WW$.
\end{Prf}

\medskip

 We expect $E'$ (resp., $E''$) to be a `maximal abelian substack' of $E$ (resp., of $E_{red}$), but we have not checked this.


 Let $\cL_{\zeta}$ be the local system on $B(\mu_n)$, the direct summand in $a_!\Qlb$, on which $\mu_n(k)$ acts by $\zeta$. Here $a\colonÊ\Spec k\to B(\mu_n)$ is the natural map. Note that $B(\mu_n)$ is a group stack, and $\cL_{\zeta}$ is a character local system on this stack. 
 
 Define a \select{$\zeta$-genuine character local system} on $E'$ as a rank one local system $\cA$ equipped with the following data. The $*$-restriction of $\cA$ to $B(\mu_n)$ is identified with $\cL_{\zeta}$. For the product map $m\colon E'\times E'\to E'$ we are given an isomorphism $\sigma\colon m^*\cA\iso \cA\boxtimes\cA$, which is associative and the restriction of $\sigma$ to $B(\mu_n)$ is compatible with the character local system structure of $\cL_{\zeta}$. 
 
 Let $a\colon E'\times E\to E$ be the product map. We have the diagram of associativity
$$
\begin{array}{ccc}
E'\times E'\times E & \toup{\id\times a} & E'\times E\\
\downarrow\lefteqn{\scriptstyle m\times\id} && \downarrow\lefteqn{\scriptstyle  a}\\
E'\times E & \toup{a} & E.
\end{array}
$$
\begin{Def} Let $(\cA,\sigma)$ be a $\zeta$-genuine character local system on $E'$. Let $Ind(\cA)$ be the category of $\Qlb$-perverse sheaves $\cF$ on $E$, on which $\mu_n(k)$ acts by $\zeta$, and which are equipped with an isomorphism
$\eta\colon a^*\cF\iso \cA\boxtimes\cF$ such that the diagram commutes
$$
\begin{array}{ccc}
(m\times\id)^*a^*\cF=(\id\times a)^*a^*\cF & \toup{(\id\times a)^*\eta} & \cA\boxtimes a^*\cF\\
\downarrow\lefteqn{\scriptstyle (m\times\id)^*\eta} && \downarrow\lefteqn{\scriptstyle \id\boxtimes\eta}\\
(m\times\id)^*(\cA\boxtimes\cF) & \toup{\sigma\boxtimes\id} & \cA\boxtimes\cA\boxtimes\cF.
\end{array}
$$
\end{Def}
 
 The group stack $E$ acts on $Ind(\cA)$ by right translations. If $n$ is odd or $\kappa(\Lambda\otimes\Lambda)\subset 2\ZZ$ then $Ind(\cA)$ is our geometric analog of the representation $\pi_x$ from Section~\ref{Section_325}. 
 
 If $n$ is even then we may repeat the construction of $\Ind(\cA)$ with $E'$ replaced by $E''$.
 
\begin{Rem} Given in addition a bilinear form $B\colon\Lambda\otimes\Lambda\to\ZZ$ with $\kappa=B+{^tB}$, one can construct (\ref{ext_T(F)_by_Gm}) as $\Gm\times T(F)$ with the group structure given by a cocycle as in Section~\ref{section_321}. In particular, the $\mu_n$-gerbe $E\to T(F)$ is trivial. 
\end{Rem}

\subsubsection{Twisted spherical sheaves}  Recall that (\ref{ext_T(F)_by_Gm}) splits over $T(\cO)$. We pick such splitting. It yields a splitting of (\ref{ext_T(F)_by_B(mu_n)}) over $T(\cO)$. So, we may consider the category $\Sph$ of perverse sheaves on $E$, which are left and right equivariant with respect to $T(\cO)$ and on which $\mu_n$ acts by $\zeta\colon \mu_n(k)\to\Qlb$. It is naturally equipped with a monoidal category structure given by the convolution. This is a geometric analog of the Hecke algebra $\cH_x$ defined in Section~\ref{section_323}. 

 The exact sequence $1\to\mu_n\to\Gm\toup{\pi_n}\Gm\to 1$, where $\pi_n(z)=z^n$, yields a morphism of group stacks $\Gm\to B(\mu_n)$, we also denote by $\cL_{\zeta}$ the restriction of $\cL_{\zeta}$ under the latter map. This is the direct summand in $(\pi_n)_!\Qlb$, on which $\mu_n(k)$ acts by $\zeta$. We may view $\Sph$ as the category of $(\Gm, \cL_{\zeta})$-equivariant perverse sheaves on $\cE$, which are also left and right $T(\cO)$-equivariant.

\begin{Rem} The monoidal category $\Sph$ has been studied in \cite{R}. One should be careful using \cite{R}, for example, (\cite{R}, Proposition~II.3.6) is wrong as stated. Besides, in \cite{R} only the case when $k$ is of characteristic zero is considered.
\end{Rem}

\smallskip

 For $\mu\in\Lambda$ let $E^{\mu}$ be the connected component of $E$ containing the preimage of $t^{\mu}\in T(F)$. It is easy to see that there is a nonzero object of $\Sph$ supported on $E^{\mu}$ iff $\mu\in\Lambda^{\sharp}$. Pick a section of (\ref{ext_T(F)_by_B(mu_n)}) over $\Lambda^{\sharp}$, that is, a morphism of group stacks $s\colon \Lambda^{\sharp}\to E$ extending the inclusion $\Lambda^{\sharp}\to T(F)$, $\mu\mapsto t^{\mu}$. The functor $s^*$ yields an equivalence of monoidal categories $\Sph\iso \Rep(\check{T}^{\sharp})$.

\subsection{Global setting}
\label{section_Global setting}

\subsubsection{} Keep the notations of Section~\ref{section_421}. According to Weissman \cite{W, W2}, the input data for the twisted Langlands correspondence for a torus should be an object $(\kappa, \tilde\Lambda)$ of $\cE(T, K_2)$ and an integer $n\ge 1$. Since $\cE(T, K_2)\subset \cE^s(T)$, it would produce an object of $\cP^{\theta}(X,\Lambda)$ by Lemma~\ref{Lm_functor_from_Weissman_to_theta_data}.

 We consider a bit more general situation. We take as initial input data an object $\theta=(\kappa,\lambda,c)\in \cP^{\theta}(X,\Lambda)$, and assume $\kappa\colon \Lambda\otimes\Lambda\to\ZZ$ to be even. So, for each $\gamma\in\Lambda$ the line bundle $\lambda^{\gamma}$ on $X$ is of parity zero. Write also $\lambda$ for the line bundle on $\Bun_T$ obtained applying the functor (\ref{functor_from_cP^theta(X,Lambda)_to_cPic(Bun_T)}) to $\theta$. Recall that $AJ^*\lambda$ is equipped with a factorizable structure. Note that $\lambda_{\cF^0}$ is trivialized. Here $\cF^0$ is the trivial $T$-torsor on $X$. The line bundle $\lambda$ on $\Bun_T$ is purely of parity zero. 
 
 As in Section~\ref{section_421}, the line bundle $\lambda$ on $\Bun_T$ yields the map $\delta_{\lambda}\colon \Lambda\to\check{\Lambda}$. By (\cite{BD2}, Lemma~4.9.2), the map $\delta_{\lambda}$ is affine with the linear part $\kappa$. We assume in addition that $\delta_{\lambda}=\kappa$. So, our $\theta=(\kappa,\lambda,c)$ is defined uniquely up to an action of the groupoid of $\check{T}$-torsors of degree zero.   
 
\subsubsection{} 
\label{Section_522}
Pick $n>1$ invertible in $k$, pick a primitive character $\zeta\colon \mu_n(k)\to\Qlb^*$. 

 Let $\Bunt_{T,\lambda}$ be the gerbe of $n$-th roots of $\lambda$ over $\Bun_T$. It classifies $\cF_T\in\Bun_T$, a $\ZZ/2\ZZ$-graded line $\cU$ of parity zero, and an isomorphism $\cU^n\iso \lambda_{\cF_T}$ of super $k$-vector spaces. For $\mu\in \Lambda$ we write $\Bunt_{T,\lambda}^{\mu}$ for the preimage of $\Bun_T^{\mu}$ in $\Bunt_{T,\lambda}$.   
Let $\D_{\zeta}(\Bunt_{T,\lambda}^{\mu})$ be the bounded derived category of $\Qlb$-sheaves on $\Bunt_{T,\lambda}^{\mu}$, on which $\mu_n(k)$ acts via $\zeta$. We have used here the natural action of $\mu_n(k)$ on $\Bunt_{T,\lambda}$ by 2-automorphisms. 
 
 Let $\Lambda^{\sharp}$ be given by (\ref{def_Lambda^sharp}), so $\Lambda^{\sharp}\subset\Lambda$ is of finite index. Set $T^{\sharp}=\Lambda^{\sharp}\otimes\Gm$. Let $i\colon T^{\sharp}\to T$ be the corresponding isogeny. Let $i_X\colon \Bun_{T^{\sharp}}\to\Bun_T$ be the push forward map. Set $\lambda^{\sharp}=i_X^*\lambda$. Write $\Bunt_{T^{\sharp},\lambda}$ for the restriction of the gerbe $\Bunt_{T,\lambda}$ under $i_X$.
 
 Let $\kappa^{\sharp}$ (resp., $\theta^{\sharp}$) denote the restriction of $\kappa$ (resp., of $\theta$) to $\Lambda^{\sharp}$. Since both $\kappa^{\sharp}$ and $\delta_{\lambda^{\sharp}}$
are divisible by $n$, we may and do pick an object $(\frac{\kappa^{\sharp}}{n}, \tau, c^{\sharp})
\in \cP^{\theta}(X,\Lambda^{\sharp})$ and an isomorphism 
\begin{equation}
\label{iso_nth_root_from_lambda}
(\frac{\kappa^{\sharp}}{n}, \tau, c^{\sharp})^n\iso \theta^{\sharp}
\end{equation}
in  $\cP^{\theta}(X,\Lambda^{\sharp})$. Note that $(\frac{\kappa^{\sharp}}{n}, \tau, c^{\sharp})$ is defined uniquely up to an action of $\check{T}^{\sharp}$-torsors on $X$, whose $n$-th power is trivialized. 

 If $n$ is odd then for any $\gamma\in\Lambda$ the line bundle $\tau^{\gamma}$ is of parity zero. If $n$ is even then it may be indeed a super line bundle. Write also $\tau$ for the super line bundle on $\Bun_{T^{\sharp}}$ obtained by applying the functor (\ref{functor_from_cP^theta(X,Lambda)_to_cPic(Bun_T)}) to $(\frac{\kappa^{\sharp}}{n}, \tau, c^{\sharp})\in \cP^{\theta}(X,\Lambda^{\sharp})$ with $\Lambda$ replaced by $\Lambda^{\sharp}$. It is equipped with a $\ZZ/2\ZZ$-graded isomorphism $\tau^n\iso \lambda^{\sharp}$ over $\Bun_{T^{\sharp}}$ obtained from (\ref{iso_nth_root_from_lambda}). This yields a section 
$$
\gs\colon \Bun_{T^{\sharp}}\to \Bunt_{T^{\sharp}, \lambda}
$$ 
of the gerbe $\Bunt_{T^{\sharp}, \lambda}\to \Bun_{T^{\sharp}}$.  

  A point of $\Bunt_{T^{\sharp},\lambda}$ is given by $\cF^{\sharp}\in\Bun_{T^{\sharp}}$ for which we set $\cF=\cF^{\sharp}\times_{T^{\sharp}} T$, and a line $\cU$ equipped with $\sigma\colon \cU^n\iso \lambda_{\cF}$. Let
$$
\pi\colon \Bunt_{T^{\sharp}, \lambda}\to \Bunt_{T,\lambda}
$$ 
be the map sending the above point to $(\cF, \cU,\sigma)$. The map $\gs$ is given by $\cU=\tau_{\cF^{\sharp}}$. 

 For $\mu\in\Lambda^{\sharp}$ we similarly define the category $\D_{\zeta}(\Bunt_{T^{\sharp}, \lambda}^{\mu})$, this is the category of objects in $\D(\Bunt_{T^{\sharp}, \lambda}^{\mu})$ on which $\mu_n(k)$ acts by $\zeta$. The section $\gs$ defines an equivalence
$$
\gs^*\colon \D_{\zeta}(\Bunt_{T^{\sharp}, \lambda}^{\mu})
\iso \D(\Bun_{T^{\sharp}}^{\mu}) \ .
$$
  
   Note that for $K\in \D_{\zeta}(\Bunt_{T^{\sharp}, \lambda}^{\mu})$ the object $\pi_!K\in \D(\Bunt_{T,\lambda}^{\mu})$ actually lies in the subcategory $\D_{\zeta}(\Bunt_{T,\lambda}^{\mu})$.
   
\begin{Lm} 
\label{Lm_some_categories_vanish}
Let $L$ be a line bundle on $B(T)$ such that $T$ acts on it via the character $\check{\lambda}\in\check{\Lambda}$. Let $\tilde B(T)$ be the gerbe of $n$-th roots of $L$ over $B(T)$, so $\tilde B(T)\to B(T)$ is a $\mu_n$-gerbe. Let $\D_{\zeta}^b(\tilde B(T))$ be the bounded derived category of $\Qlb$-complexes on $\tilde B(T)$ on which $\mu_n(k)$ acts by $\zeta\colon \mu_n(k)\to\Qlb^*$. If $\check{\lambda}\notin n\check{\Lambda}$ then $\D_{\zeta}^b(\tilde B(T))=0$.
\end{Lm}
\begin{Prf}
Let $G$ be the kernel of the map $T\times \Gm\to \Gm$, $(t,z)\mapsto \check{\lambda}(t)z^{-n}$. Then $\tilde B(T)\iso B(G)$ naturally. Consider the weight lattice $\check{\Lambda}\oplus\ZZ$ of $T\times\Gm$. The above map is the weight $\check{\lambda}-n\check{\epsilon}$, where $\check{\epsilon}$ is the standard weight of $\Gm$. Let $m\ge 1$ be the biggest integer such that $\frac{1}{m}(\check{\lambda}-n\check{\epsilon})\in \check{\Lambda}\oplus\ZZ$. This is a divisor of $n$. In particular $m$ is invertible in $k$. By Remark~\ref{Rem_pi_0_of_Ker_of_weight} below, $\pi_0(G)\iso \mu_m(k)$. By our assumption, $m<n$. The action of $\mu_n(k)$ by 2-automorphisms on $B(G)$ factors through the map $\mu_n(k)\to \pi_0(G)$, so $\mu_n(k)$ can not act via a faithfull character.
\end{Prf}

\begin{Rem} 
\label{Rem_pi_0_of_Ker_of_weight}
If $\check{\lambda}\in\check{\Lambda}$, consider the kernel $\Ker\check{\lambda}$ of $\check{\lambda}\colon T\to\Gm$. Let $\check{\mu}\in\check{\Lambda}$ be such that $\check{\Lambda}/\ZZ\check{\mu}$ is torsion free, and $\check{\lambda}=m\check{\mu}$ for some $m\ge 1$. The group $\pi_0(\Ker\check{\lambda})$ of connected components of $\Ker\check{\lambda}$ is $\pi_0(\mu_m(k))$. If $m$ is invertible in $k$ then $\pi_0(\mu_m(k))=\mu_m(k)$. 
\end{Rem}
   
\begin{Prf}\select{of Proposition~\ref{Pp_intro_vanishing}} \  \  Pick a $k$-point $\cF\in \Bun_T^{\mu}$. It gives a map $B(T)\to \Bun_T^{\mu}$, as $T$ is the group of automorphisms of $\cF$. Let $L$ be the restriction of the line bundle $\lambda$ to $B(T)$. Let $\tilde B(T)$ be the gerbe of $n$-th roots of $L$ over $B(T)$. It suffices to show that $\D_{\zeta}^b(\tilde B(T))=0$. The group $T$ acts on $L$ by the character $\delta_{\lambda}(\mu)=\kappa(\mu)\in \check{\Lambda}$. The condition $\mu\notin \Lambda^{\sharp}$ is equivalent to $\kappa(\mu)\notin n\check{\Lambda}$. So, $\D_{\zeta}^b(\tilde B(T))=0$ by Lemma~\ref{Lm_some_categories_vanish}. 
\end{Prf}

\subsubsection{Hecke functors} 
\label{section_Hecke_functors}
In this section we construct an action of $\Bun_{T^{\sharp}}$ on $\Bunt_{T,\lambda}$.

 A point of $\Bunt_{T,\lambda}$ is a pair $(\cF,\cU)$, where $\cF\in\Bun_T$, $\cU$ is a line of parity zero, and $\cU^n\iso \lambda_{\cF}$ is a $\ZZ/2\ZZ$-graded isomorphism. For $\mu\in\Lambda^{\sharp}$ we define a map $m_{\mu}\colon X\times \Bunt_{T,\lambda}\to\Bunt_{T,\lambda}$ as follows. It sends $(x\in X, (\cF,\cU)\in\Bunt_{T,\lambda})$ to $(\cF(\mu x), \cU')$, where 
\begin{equation}
\label{def_of_cB'_for_m_mu}
\cU'=\cU\otimes  (\cL^{\frac{\kappa(\mu)}{n}}_{\cF})_x\otimes \tau_{\cO(\mu x)}
\end{equation}
is equipped with the isomorphism 
\begin{equation}
\label{iso_for_cB'}
(\cU')^n\iso \lambda_{\cF(\mu x)}
\end{equation}
given by Proposition~\ref{Pp_factoriztion_for_twists_by_divisors}. 
Let us explain that for $\mu\in \Lambda^{\sharp}$ we may view $\cO(\mu x)$ as a $k$-point of $\Bun_{T^{\sharp}}$, and $\tau_{\cO(\mu x)}$ denotes here the fibre of $\tau$ at $\cO(\mu x)\in \Bun_{T^{\sharp}}$. Recall that over $\Bun_{T^{\sharp}}$ we have an isomorphism $\tau^n\iso \lambda\mid_{\Bun_{T^{\sharp}}}$ obtained from (\ref{iso_nth_root_from_lambda}). This is how (\ref{iso_for_cB'}) is obtained. 

 The Hecke functor 
$$
\H^{\mu}\colon \D_{\zeta}(\Bunt_{T, \lambda})\to \D_{\zeta}(X\times \Bunt_{T, \lambda})
$$ 
is defined as $\H^{\mu}(K)=m_{\mu}^*K$. For $a\in \mu_n(k)$ the corresponding 2-automorphism of $X\times\Bunt_{T,\lambda}$ acts by $a$ on $\cU$ and trivially on $(\cF, x)$. The image of this 2-automorphism under $m_{\mu}$ acts by $a$ on $\cU'$ and trivially on $\cF(\mu x)$. So, if $K\in \D_{\zeta}(\Bunt_{T, \lambda})$ then $\mu_n(k)$ acts on $m_{\mu}^*K$ by $\zeta$.

  If $\Lambda'$ is a free abelian group of finite type and $j\colon \Lambda\to\Lambda'$ is a linear map, for $T'=\Lambda'\otimes\Gm$ we get a map $\bar  j\colon \Bun_T\to\Bun_{T'}$ such that if $\check{\lambda}\in\check{\Lambda}'$ then $\cL^{\check{\lambda}}_{\bar j(\cF)}\iso \cL^{\check{\lambda}\comp j}_{\cF}$. Here $\check{\Lambda}'=\Hom(\Lambda',\ZZ)$. The diagram 
$$
\begin{array}{ccc}
\Lambda^{\sharp} & \hook{} & \Lambda\\
\downarrow\lefteqn{\scriptstyle \frac{\kappa}{n}} && \downarrow\lefteqn{\scriptstyle\kappa}\\
\check{\Lambda} & \toup{n} & \check{\Lambda}
\end{array}
$$
yields a diagram of morphisms 
$$
\begin{array}{ccc}
\Bun_{T^{\sharp}} & \toup{i_X} & \Bun_T\\
\downarrow\lefteqn{\scriptstyle i_{\kappa} } && \downarrow\lefteqn{\scriptstyle  \bar\kappa}\\
\Bun_{\check{T}} & \toup{n_X} & \Bun_{\check{T}},
\end{array}
$$
where $n_X$ sends $\cF$ to $\cF^{\otimes n}$.  
\begin{Lm} 
\label{Cor_action_finally}
For $\cF\in\Bun_T, \cT^{\sharp}\in\Bun_{T^{\sharp}}$ let $\cT=\cT^{\sharp}\times_{T^{\sharp}} T$. Our choice of (\ref{iso_nth_root_from_lambda}) yields an isomorphism
$$
\lambda_{\cF\otimes \cT}\iso \lambda_{\cF}\otimes \tau^{\otimes n}_{\cT^{\sharp}}\otimes (\cL^{univ}_{i_{\kappa}(\cT^{\sharp}), \cF})^{\otimes n} \ .
$$
\end{Lm}
\begin{Prf}
By Proposition~\ref{Pp7_factorization_of_lambda_on_BunT}, $\lambda_{\cF\otimes \cT}\iso \lambda_{\cF}\otimes \tau^n_{\cT^{\sharp}}\otimes \cL^{univ}_{\bar\kappa(\cT),\cF}$. Now $\cL^{univ}_{\bar\kappa(\cT),\cF}\iso \cL^{univ}_{i_{\kappa}(\cT^{\sharp})^{\otimes n}, \cF}\iso (\cL^{univ}_{i_{\kappa}(\cT^{\sharp}), \cF})^{\otimes n}$ because of the biextension structure on $\cL^{univ}$. 
\end{Prf}

\medskip

 Let $a\colon \Bun_{T^{\sharp}}\times\Bunt_{T,\lambda}\to \Bunt_{T,\lambda}$ be the map sending $(\cT^{\sharp}, \cF, \cU, \cU^n\iso \lambda_{\cF})$ to 
$
(\cF\otimes\cT, \cU'),
$
where $\cT=\cT^{\sharp}\times_{T^{\sharp}} T$, and 
\begin{equation}
\label{def_cU'_for_the_map_a}
\cU'=\cU\otimes \tau_{\cT^{\sharp}}\otimes \cL^{univ}_{i_{\kappa}(\cT^{\sharp}),\cF}
\end{equation}
is equipped with the isomorphism $\cU'^n\iso \lambda_{\cF\otimes\cT}$ given by Lemma~\ref{Cor_action_finally}. 

\begin{Lm} The map $a$ has a natural structure of an action of the group stack $\Bun_{T^{\sharp}}$ on $\Bunt_{T,\lambda}$.
\end{Lm}
\begin{Prf}
The bilinear form associated to $\tau$ is $\frac{\kappa^{\sharp}}{n}\colon \Lambda^{\sharp}\otimes\Lambda^{\sharp}\to\ZZ$. Using Remark~\ref{Rem_some_symmetry} and applying Proposition~\ref{Pp7_factorization_of_lambda_on_BunT} to $\tau$ on $\Bun_{T^{\sharp}}$, for $\cT^{\sharp}, \cG^{\sharp}\in\Bun_{T^{\sharp}}$ one gets an isomorphism 
\begin{equation}
\label{tau_factorization}
\tau_{\cT^{\sharp}}\otimes \tau_{\cG^{\sharp}}\otimes\cL^{univ}_{i_{\kappa}(\cG^{\sharp}), \cT}\iso \tau_{\cG^{\sharp}\otimes\cT^{\sharp}},
\end{equation}
where $\cT=\cT^{\sharp}\times_{T^{\sharp}} T$. This combined with Lemma~\ref{Cor_action_finally} gives a 2-morphism making the following diagram 2-commutative
$$
\begin{array}{ccc}
\Bun_{T^{\sharp}}\times\Bun_{T^{\sharp}}\times\Bunt_{T,\lambda} & \toup{\id\times a} &\Bun_{T^{\sharp}}\times\Bunt_{T,\lambda} \\
\downarrow{\lefteqn{\scriptstyle m\times \id}} && \downarrow{\lefteqn{\scriptstyle a}}\\
\Bun_{T^{\sharp}}\times\Bunt_{T,\lambda} & \toup{a} &\Bunt_{T,\lambda} .
\end{array}
$$
Here $m$ is the product map for $\Bun_{T^{\sharp}}$. Besides, there is a natural 2-morphism identifying the restriction of $a$ to the trivial $T^{\sharp}$-torsor with the identity map.
\end{Prf}
 
\medskip

 Consider the map $e^0\colon B(\mu_n)\to \Bunt_{T,\lambda}$ sending $(\cU, \cU^n\iso k)$ to $(\cF^0, \cU, \cU^n\iso \lambda_{\cF^0})$. We used the fact that $\lambda_{\cF^0}\iso k$ canonically, here $\cF^0$ is the trivial $T$-torsor on $X$. The composition 
$$
\Bun_{T^{\sharp}}\times B(\mu_n)\toup{\id\times e^0} \Bun_{T^{\sharp}}\times\Bunt_{T,\lambda}\toup{a}\Bunt_{T,\lambda}
$$ 
is naturally 2-isomorphic to $\pi\colon \Bunt_{T^{\sharp},\lambda}\to\Bunt_{T,\lambda}$. We have a 2-commutative diagram
\begin{equation}
\label{diag_for_a}
\begin{array}{ccc} 
\Bun_{T^{\sharp}}\times\Bunt_{T,\lambda}& \toup{\bar a} & \Bun_{T^{\sharp}}\times\Bunt_{T,\lambda}\\
& \searrow\lefteqn{\scriptstyle a} & \downarrow\lefteqn{\scriptstyle{\pr_2}}\\
&&\Bunt_{T,\lambda},
\end{array}
\end{equation}
where $\bar a$ is the isomorphism sending $(\cT^{\sharp}, \cF,\cU, \cU^n\iso \lambda_{\cF})$ to $(\cT^{\sharp}, \cF\otimes\cT, \cU', \cU'^n\iso \lambda_{\cF\otimes\cT})$ with $\cT=\cT^{\sharp}\times_{T^{\sharp}} T$, and $\cU'$ is given by (\ref{def_cU'_for_the_map_a}). So, $a$ is a $\Bun_{T^{\sharp}}$-torsor.

   Let $\check{T}^{\sharp}$ be the Langlands dual to $T^{\sharp}$ torus  over $\Qlb$, let $E$ be a $\check{T}^{\sharp}$-local system on $X$. Let $AE$ denote the corresponding automorphic local system on $\Bun_{T^{\sharp}}$. For a $\Lambda^{\sharp}$-valued divisor $D=\sum_{x\in X} \,\mu_x x$ on $X$, write $\cO(D):=AJ(D)\in\Bun_{T^{\sharp}}$ then
$$
AE_{\cO(D)}\iso \bigotimes_{x\in X} E^{\mu_x}_x
$$ 
canonically. Here for $\mu\in\Lambda^{\sharp}$ we denoted by $E^{\mu}$ the push forward of $E$ via $\mu\colon \check{T}^{\sharp}\to \Gm$. Note that $AE$ is a character local system on $\Bun_{T^{\sharp}}$. So, for the product map $m$ there is a natural isomorphism $m^*AE\iso  AE\boxtimes AE$, in partucular, the restriction of $AE$ to the trivial $T^{\sharp}$-torsor is trivialized.

\begin{Def} A $E$-Hecke eigensheaf in $\D_{\zeta}(\Bunt_{T,\lambda})$ is an object $K\in \D_{\zeta}(\Bunt_{T,\lambda})$ equipped with a $(\Bun_{T^{\sharp}}, AE)$-equivariant structure. This means that it is equipped with an isomorphism $a^*K\iso AE\boxtimes K$, which is associative, and whose restriction to the unit section is trivialized (in a way compatible with the character local system structure on $AE$). 
\end{Def}

 Note that if $K\in \D_{\zeta}(\Bunt_{T,\lambda})$ is a $E$-Hecke eigen-sheaf then for $\mu\in\Lambda^{\sharp}$ we get an isomorphism $\H^{\mu}(K)\iso E^{\mu}\boxtimes K$. 
 
 Let 
$$
a^{\sharp}\colon \Bun_{T^{\sharp}}\times\Bunt_{T^{\sharp},\lambda}\to\Bunt_{T^{\sharp},\lambda}
$$ 
be the map sending $(\cT^{\sharp}, \cF^{\sharp}, \cU, \cU^n\iso \lambda_{\cF})$ with $\cT=\cT^{\sharp}\times_{T^{\sharp}}T$, $\cF=\cF^{\sharp}\times_{T^{\sharp}}T$ to the collection $(\cT^{\sharp}\otimes\cF^{\sharp}, \cU')$, where $\cU'$ is given by (\ref{def_cU'_for_the_map_a}) and equipped with the isomorphism $\cU'^n\iso \lambda_{\cF\otimes\cT}$ given by  
Lemma~\ref{Cor_action_finally}. The map $a^{\sharp}$ is equipped with a structure of an action map of $\Bun_{T^{\sharp}}$ on $\Bunt_{T^{\sharp},\lambda}$. The diagram is naturally 2-commutative and cartesian
\begin{equation}
\label{diag_cartesian_for_Hecke}
\begin{array}{ccc}
\Bun_{T^{\sharp}}\times\Bunt_{T^{\sharp},\lambda}& \toup{a^{\sharp}} &\Bunt_{T^{\sharp},\lambda}\\
\downarrow\lefteqn{\scriptstyle \id\times\pi} && \downarrow\lefteqn{\scriptstyle\pi}\\
\Bun_{T^{\sharp}}\times\Bunt_{T,\lambda}& \toup{a}&\Bunt_{T,\lambda} .
\end{array}
\end{equation}

 Define the derived category $\D_{\zeta}(\Bun_{T^{\sharp}}\times\Bunt_{T^{\sharp},\lambda})$ similarly, that is, by requiring that $\mu_n(k)$ acts by $\zeta$. Consider the 2-automorphism of $\Bun_{T^{\sharp}}\times\Bunt_{T^{\sharp},\lambda}$ acting by $a\in \mu_n(k)$ on $\cU$ and trivially on $\cT^{\sharp}, \cF^{\sharp}$. Its image under $a^{\sharp}$ acts by $a$ on $\cU'$ and trivially on $\cT^{\sharp}\otimes\cF^{\sharp}$. Therefore, $(a^{\sharp})^*\colon \D_{\zeta}(\Bunt_{T^{\sharp},\lambda})\to \D_{\zeta}(\Bun_{T^{\sharp}}\times\Bunt_{T^{\sharp},\lambda})$. 
 
\medskip\noindent  
\begin{Prf}\select{of Proposition~\ref{Pp_H_eigensheaves} (i)} \ \ 
The local system $W$ is naturally equipped with an isomorphism $(a^{\sharp})^*W\iso  AE\boxtimes W$, which is associative, and its restriction to the unit section of $\Bun_{T^{\sharp}}$ is trivialized. Since (\ref{diag_cartesian_for_Hecke}) is cartesian, $\pi_!W$ is a $E$-Hecke eigen-sheaf naturally.
\end{Prf}

\subsubsection{} 
\label{Section_524}
As in Section~\ref{section_323}, define $K$ by the exact sequence $1\to K\to T^{\sharp}\toup{i}T\to 1$.
\begin{Lm} 
\label{Lm_K-torsors_equivariance}
Let $\cG_1$ be a $K$-torsor on $X$, $\mu\in\Lambda^{\sharp}$. The $\mu_n$-torsor on $\Bun^{\mu}_{T^{\sharp}}$ with fibre $\tau_{\cG\otimes \cG_1}\otimes\tau_{\cG}^{-1}$ at $\cG\in \Bun^{\mu}_{T^{\sharp}}$ is constant on $\Bun^{\mu}_{T^{\sharp}}$ and independent of $\mu$. We have canonically $\tau_{\cG\otimes \cG_1}\otimes\tau_{\cG}^{-1}\iso \tau_{\cG_1}$. 
\end{Lm}
\begin{Prf}
The $n$-th power of the line bundle on $\Bun_{T^{\sharp}}^{\mu}$ with fibre $\tau_{\cG\otimes \cG_1}\otimes\tau_{\cG}^{-1}$ at $\cG$ is trivialized, so we think of it as a $\mu_n$-torsor. By Proposition~\ref{Pp7_factorization_of_lambda_on_BunT}, for $\cG\in\Bun_{T^{\sharp}}$ we get $\tau_{\cG\otimes\cG_1}\iso \tau_{\cG}\otimes\tau_{\cG_1}\otimes (^{\frac{\kappa^{\sharp}}{n}}\cL^{univ}_{\cG,\cG_1})$. The map $\frac{\kappa^{\sharp}}{n}\colon \Lambda^{\sharp}\times\Lambda^{\sharp}\to\ZZ$ factors as $\Lambda^{\sharp}\times\Lambda^{\sharp}\toup{\id\times i}\Lambda^{\sharp}\times\Lambda\to \ZZ$. Therefore, 
$$
^{\frac{\kappa^{\sharp}}{n}}\cL^{univ}_{\cG,\cG_1}\iso {^{\frac{\kappa}{n}}\cL^{univ}_{\cG, \cF^0}}\iso k,
$$
where $\cF^0$ is the trivial $T$-torsor. The latter isomorphism comes from the biextension structure on $\cL^{univ}$. 
\end{Prf}

\medskip

 Write $\und{\Bun}_{T^{\sharp}}^{\mu}$ for the coarse moduli space of $\Bun_{T^{\sharp}}^{\mu}$, similarly for $\und{\Bun}_T^{\mu}$. For $\mu\in \Lambda^{\sharp}$ and a $k$-point $\eta\in \und{\Bun}_{T^{\sharp}}^{\mu}$, $\pi_1(\eta, \und{\Bun}_{T^{\sharp}}^{\mu})$ is abelian and independent of $\eta$ and $\mu$ up to a canonical isomorphism, we simply denote it by $\pi_1(\und{\Bun}_{T^{\sharp}})$. Similarly for $\pi_1(\und{\Bun}_T)$. The map $i_X$ induces a Galois covering $\und{\pi}\colon \und{\Bun}_{T^{\sharp}}^{\mu}\to \und{\Bun}_T^{\mu}$ with Galois group $\H^1(X, K)$. It yields an exact sequence
\begin{equation}
\label{ext_H1(X,K)_by_pi_1}
1\to \pi_1(\und{\Bun}_{T^{\sharp}})\to \pi_1(\und{\Bun}_T)\to \H^1(X, K)\to 1.
\end{equation}
 
  Note that over $\Bun_{T^{\sharp}}^0$ the line bundle $\tau$ descends to a line bundle on its coarse moduli space $\uBun_{T^{\sharp}}^0$, which is an abelian variety. Let $K(\tau)$ denote the kernel of the map $\phi_{\tau}\colon A\to \hat A$ defined as in (\cite{P}, Corollary~8.6), here $A=\uBun_{T^{\sharp}}^0$. By Lemma~\ref{Lm_K-torsors_equivariance}, $\H^1(X, K)\subset K(\tau)$. This inclusion may be strict, see example in Section~\ref{section_case_T=Gm} below.

 If $S$ is a scheme, $H$ is a flat finitely presented and separated group scheme over $S$, $\cX$ is an algebraic stack over $S$, and $\alpha\colon \pr_2\to\pr_2$ is a 2-morphism for the projection $\pr_2\colon H\times_S\cX\to\cX$ defining an action of $H$ on $\cX$ by 2-automorphisms, we will use a rigidification of $\cX$ along $H$ (\cite{ACV}, Definition~5.1.9 and Theorem~5.1.5) obtained by killing $H$ inside the automorphisms of objects of $\cX$.

  Let $'\Bun_{T^{\sharp}}$ be obtained from $\Bun_{T^{\sharp}}$ by
killing the group $\H^0(X, K)$ inside the automorphisms of $\Bun_{T^{\sharp}}$. The projection $\Bun_{T^{\sharp}}\to {'\Bun_{T^{\sharp}}}$ is a $K$-gerbe. One has non-canonically $'\Bun_{T^{\sharp}}\iso \uBun_{T^{\sharp}}\times B(T)$. The map $i_X$ factors through $'i_X\colon {'\Bun_{T^{\sharp}}}\to \Bun_T$, the map $'i_X$ is a Galois covering with Galois group $\H^1(X,K)$. 

  Since $\delta_{\tau}=\frac{\kappa^{\sharp}}{n}$, it follows that $K$ acts trivially on $\tau$. So, $\tau$ descends to a line bundle $'\tau$ on $'\Bun_{T^{\sharp}}$. The isomorphism $\tau^n\iso \lambda^{\sharp}$ also descends to an isomorphism 
\begin{equation} 
\label{iso_over'Bun} 
'\tau^n\iso ('i_X)^*\lambda
\end{equation}
on $'\Bun_{T^{\sharp}}$. Let $'\pi\colon {'\Bun_{T^{\sharp}}}\to \Bunt_{T,\lambda}$ be the map corresponding to (\ref{iso_over'Bun}). Now applying Remark~\ref{Rem_for_Section_B1} ii) and Lemma~\ref{Lm_K-torsors_equivariance} to the Galois covering $'i_X\colon {'\Bun_{T^{\sharp}}}\to \Bun_T$, we get a central extension
\begin{equation}
\label{ext_H^1(X,K)_by_mu_n}
1\to \mu_n(k)\to \Gamma\to \H^1(X, K)\to 1
\end{equation}
(eventually depending on $\mu$) such that $'\pi\colon {'\Bun_{T^{\sharp}}^{\mu}}\to
\Bunt_{T,\lambda}^{\mu}$
is a Galois covering with Galois group $\Gamma$, here $\Gamma$ acts on $'\Bun_{T^{\sharp}}$ via its quotient $\H^1(X,K)$. 
We write $\Gamma_{\mu}=\Gamma$ if we need to express the dependence of $\Gamma$ on $\mu$. 

  Given $\nu\in\Lambda^{\sharp}$ and a $k$-point $\cT^{\sharp}\in\Bun_{T^{\sharp}}^{\nu}$, the action by $\cT^{\sharp}$ yields a diagram
$$
\begin{array}{ccc}
'\Bun^{\mu}_{T^{\sharp}} & \to & {'\Bun^{\mu+\nu}_{T^{\sharp}}}\\
\downarrow\lefteqn{\scriptstyle '\pi} && \downarrow\lefteqn{\scriptstyle '\pi}\\
\Bunt^{\mu}_{T,\lambda} & \to & \Bunt^{\mu+\nu}_{T,\lambda},
\end{array}
$$  
where the horizontal arrows are isomorphisms. This provides an isomorphism $\alpha_{\nu}\colon \Gamma_{\mu}\iso \Gamma_{\mu+\nu}$, which depends only on $\nu,\mu$ and not on $\cT^{\sharp}$, because $\Bun_{T^{\sharp}}^{\nu}$ is connected. Moreover, for $\nu_i\in\Lambda^{\sharp}$ we get $\alpha_{\nu_1}\alpha_{\nu_2}=\alpha_{\nu_1+\nu_2}$. In this sense $\Gamma_{\mu}$ is independent of $\mu$. 
  
  Define the category of  \select{$\zeta$-genuine} local systems on $\Bunt_{T,\lambda}$ as the category of local systems in $\D_{\zeta}(\Bunt_{T,\lambda})$. 
For $\mu\in\Lambda^{\sharp}$ let 
$$
1\to \mu_n(k)\to \bar\Gamma\to \pi_1(\und{\Bun}_T)\to 1
$$
be the exact sequence obtained as the pull-back of (\ref{ext_H^1(X,K)_by_mu_n}) by (\ref{ext_H1(X,K)_by_pi_1}). 
  
\begin{Cor} 
\label{Cor_first}
Let $\mu\in\Lambda^{\sharp}$. \\
(1) The category of $\zeta$-genuine $\Qlb$-local systems on $\Bunt_{T,\lambda}^{\mu}$ is equivalent to the category of $\Gamma_{\mu}$-equivariant $\Qlb$-local systems $V$ on $'\Bun_{T^{\sharp}}^{\mu}$ such that the subgroup $\mu_n(k)$ acts on $V$ by the character $\zeta\colon \mu_n(k)\to\Qlb^*$. \\
(2) The category of $\zeta$-genuine $\Qlb$-local systems on $\Bunt_{T,\lambda}^{\mu}$ is equivalent to the category of finite-dimensional representations of $\bar\Gamma$, on which $\mu_n(k)$ acts by $\zeta$.
\end{Cor}  
\begin{Prf}
(2) Let $\Gamma$ act on $\und{\Bun}^{\mu}_{T^{\sharp}}$ via its quotient $\H^1(X,K)$. The natural map $\gr\colon {'\Bun^{\mu}_{T^{\sharp}}}\to \und{\Bun}^{\mu}_{T^{\sharp}}$ is $\Gamma$-equivariant. Now $\gr^*$ gives an equivalence of the category of $\Gamma$-equivariant $\Qlb$-local system on $\und{\Bun}^{\mu}_{T^{\sharp}}$, on which $\mu_n(k)$ acts by $\zeta$ with the category of $\Gamma$-equivariant $\Qlb$-local system on $'\Bun^{\mu}_{T^{\sharp}}$, on which $\mu_n(k)$ acts by $\zeta$. Indeed, for any $\Gm$-gerbe $\alpha: \tilde Y\to Y$ over a stack $Y$, $\alpha^*$ is an equivalence between the categories of $\Qlb$-local systems on $Y$ and on $\tilde Y$. 
\end{Prf}

\medskip

    In (\cite{P}, Section~10.4) the Weil pairing on $K(\tau)$ has been constructed, it is a skew-symmetric bilinear form $b_{\tau}\colon K(\tau)\times K(\tau)\to\Gm$ associated to the above map $\phi_{\tau}\colon A\to \hat A$ with $A=\uBun_{T^{\sharp}}^0$. By our construction, the commutator of the extension (\ref{ext_H^1(X,K)_by_mu_n}) is the restriction of $b_{\tau}$ to $\H^1(X, K)$. The latter takes values in $\mu_n$, because $\tau^n$ descends under $\Bun_{T^{\sharp}}\to\Bun_T$. 
    
    The map $\phi_{\tau}$ is not always an isogeny, and $K(\tau)$ is not always a finite group scheme. If $\phi_{\tau}$ is an isogeny then by (\cite{P}, Theorem~10.1) the bilinear form $b_{\tau}$ is non-degenerate, it identifies $K(\tau)$ with its Cartier dual.

\begin{Pp} 
\label{Pp_commutator_is_nondeg}
The commutator $(\cdot,\cdot)_c\colon \H^1(X, K)\times \H^1(X, K)\to\mu_n(k)$ of (\ref{ext_H^1(X,K)_by_mu_n}) is non-degenerate and independent of $\mu\in\Lambda^{\sharp}$, it induces an isomorphism 
$$
\H^1(X, K)\iso \Hom(\H^1(X, K), \mu_n(k)) \ .
$$
\end{Pp}
\begin{Prf} The commutator $(\cdot,\cdot)_c$ is described as follows. Since $\Lambda/\Lambda^{\sharp}$ is a $\ZZ/n\ZZ$-module, we get $\Hom(\check{\Lambda}^{\sharp}/\check{\Lambda}, \frac{1}{n}\ZZ/\ZZ)\iso \Lambda/\Lambda^{\sharp}$ canonically. So, $K\iso \Hom(\check{\Lambda}^{\sharp}/\check{\Lambda}, \mu_n)\iso (\Lambda/\Lambda^{\sharp})\otimes\mu_n$ canonically. The cup-product gives the pairing 
\begin{equation}
\label{pairing_first}
\H^1(X, \Lambda\otimes\mu_n)\otimes\H^1(X, \Lambda\otimes\mu_n)\to\H^2(X, \Lambda\otimes\Lambda\otimes \mu_n^{\otimes 2})\toup{\kappa}\H^2(X, \mu_n^{\otimes 2})\iso \mu_n.
\end{equation}
We have an exact sequence $0\to (\Lambda^{\sharp}/n\Lambda)\otimes\mu_n\to \Lambda\otimes\mu_n\to (\Lambda/\Lambda^{\sharp})\otimes \mu_n\to 0$. 
It yields a long exact sequence 
$$
\begin{array}{cccc}
\H^1(X, \Lambda\otimes\mu_n)\toup{\nu} \H^1(X, K)\to & \H^2(X, (\Lambda^{\sharp}/n\Lambda)\otimes\mu_n) & \to & \H^2(X, \Lambda\otimes\mu_n)\\
& \downarrow && \downarrow\\
& (\Lambda^{\sharp}/n\Lambda) &  & \Lambda/n\Lambda,
\end{array}
$$ 
where the vertical arrows are canonical isomorphisms. This shows that $\nu$ is surjective. 
The pairing (\ref{pairing_first}) is skew-symmetric, because $\kappa$ is symmetric, and the cup-product is skew-symmetric. This pairing vanishes  
$(\Ker\nu)\times \H^1(X, \Lambda\otimes\mu_n)$, because it vanishes after restriction to $\H^1(X, (\Lambda^{\sharp}/n\Lambda)\otimes\mu_n)\otimes \H^1(X, \Lambda\otimes\mu_n)$. So, we get a non-degenerate pairing $\H^1(X, K)\times \H^1(X, K)\to\mu_n$.
\end{Prf}

\medskip\noindent
\select{Second proof of Proposition~\ref{Pp_commutator_is_nondeg} in the case $\dim T=1$}. \   
\Step 1
Assume first $\dim T$ arbitrary. Let us show that the group $K(\tau)$ identifies with the group of $\cT^{\sharp}\in \uBun_{T^{\sharp}}^0$ such that for all $\mu\in\Lambda^{\sharp}$ the line bundle $\cL^{\kappa(\mu)/n}_{\cT^{\sharp}}$ is trivial on $X$. 

 Let $\cT^{\sharp}\in \uBun_{T^{\sharp}}^0$. Consider the line bundle on $\uBun_{T^{\sharp}}^0$ with fibre $\tau_{\cT^{\sharp}\otimes \cG^{\sharp}}\otimes \tau^{-1}_{\cG^{\sharp}}$ at $\cG^{\sharp}\in \uBun_{T^{\sharp}}^0$. From (\ref{tau_factorization}) we see that this line bundle is constant if and only if the line bundle with fibre $\cL^{univ}_{i_{\kappa}(\cG^{\sharp}), \cT}$ at $\cG^{\sharp}\in \uBun_{T^{\sharp}}^0$ is constant. Here $\cT=\cT^{\sharp}\times_{T^{\sharp}} T$. This is equivalent to requiring that for all $\mu\in\Lambda^{\sharp}$ the line bundle on $X$ with fibre 
$$
\cL^{univ}_{i_{\kappa}(\cO(\mu x)), \cT}=(\cL^{\frac{\kappa(\mu)}{n}}_{\cT})_x
$$ 
at $x\in X$ is trivial. 
 
 Let $\check{\Theta}$ denote the image of $\frac{1}{n}\kappa\colon \Lambda^{\sharp}\to\check{\Lambda}$. Let $T_{\Theta}$ be the split torus, whose weights are $\check{\Theta}$. The inclusion $\check{\Theta}\subset \check{\Lambda}$ gives a morphism $T\to T_{\Theta}$. Let $K_1$ be the kernel of $T\to T_{\Theta}$, so $\H^1(X, K_1)$ is the kernel of 
$\uBun_T^0\to \uBun_{T_{\Theta}}^0$. Then $K(\tau)$ is the preimage of $\H^1(X, K_1)$ under $\uBun_{T^{\sharp}}^0\to \uBun_T^0$.  

\smallskip

\Step 2 If $\rk(\Lambda)=1$ then the only nontrivial case is $\kappa\ne 0$. In this case let $d>0$ be such that $\kappa(\Lambda\otimes\Lambda)=d\ZZ$, recall that $d$ is even. Let $e$ be the smallest positive integer such that $de\in n\ZZ$. Then $\Lambda^{\sharp}=e\Lambda$, $\check{\Theta}=\frac{ed}{n}\check{\Lambda}$.
So, $K\iso \mu_e$ and $K_1\iso \mu_{\frac{ed}{n}}$. Recall that the Weil pairing $b_{\tau}$ on $K(\tau)$ is non-degenerate. Since $e$ and $\frac{ed}{n}$ are relatively prime, the orders of the groups $\H^1(X,K)$ and $\H^1(X, K_1)$ are relatively prime. So, the restriction of $b_{\tau}$ to $\H^1(X,K)$ is also non-degenerate.
\QED 
  
\begin{Cor} The center of $\Gamma$ is $\mu_n(k)$. There is a unique irreducible representation of $\Gamma$ with central character $\zeta$. \QED
\end{Cor}

Let $\H\subset \H^1(X,K)$ be a maximal subgroup isotropic with respect to $(\cdot,\cdot)_c$. The order of $\H^1(X,K)$ is $e^{2g}$, where $e$ is the order of $\Lambda/\Lambda^{\sharp}$, so $H$ is of order $e^g$. Recall that $\H^1(X,K)$ acts on $'\Bun_{T^{\sharp}}$. Let $\Bun_{T^{\sharp}, H}$ denote the quotient of $'\Bun_{T^{\sharp}}$ by $H$. For $\mu\in \Lambda^{\sharp}$ write $\Bun_{T^{\sharp}, H}^{\mu}$ for the image of $'\Bun_{T^{\sharp}}^{\mu}$ in $\Bun_{T^{\sharp}, H}$. Let
$$
'\Bunt_{T^{\sharp}}\,\to \,\Bunt_{T^{\sharp}, H}\;\toup{\pi_H}\; \Bunt_{T,\lambda}
$$
be obtained from $'\Bun_{T^{\sharp}}\to \Bun_{T^{\sharp}, H}\to \Bun_T$ by the base change $\Bunt_{T,\lambda}\to\Bun_T$.

\medskip 

\begin{Prf}\select{of Proposition~\ref{Pp_H_eigensheaves} (ii)}  Let $\mu\in\Lambda^{\sharp}$. Let $H_{\Gamma}$ be the preimage of $H$ in $\Gamma$, $\bar H$ be the preimage of $H_{\Gamma}$ in $\bar\Gamma$. Let $\chi\colon \pi_1(\und{\Bun}_{T^{\sharp}})\to\Qlb^*$ be the character corresponding to $AE$. By (\cite{W}, Proposition~2.1), we may and do pick a character $\bar\chi\colon \bar H\to\Qlb^*$ extending 
$$
\zeta\boxtimes\chi\colon \mu_n(k)\times \pi_1(\und{\Bun}_{T^{\sharp}})\to\Qlb^*\ .
$$ 
It yields a $H_{\Gamma}$-equivariant structure on $AE$ over $'\Bun_{T^{\sharp}}^{\mu}$ such that $\mu_n(k)$ acts on it by $\zeta$. 

 Since $'\Bun_{T^{\sharp}}\to \Bunt_{T^{\sharp}, H}$ is a Galois covering with Galois group $H_{\Gamma}$, we get a rank one local system $AE_H$ on $ \Bunt_{T^{\sharp}, H}$ equipped with an isomorphism of its restriction to $'\Bun_{T^{\sharp}}$ with $AE$. Set $\cK_E=\pi_{H !} (AE_H)$.
 
 Over $\Bunt_{T,\lambda}^{\mu}$ the local system $\pi_!W$ of Proposition~\ref{Pp_H_eigensheaves} corresponds via the equivalence of Corollary~\ref{Cor_first} to the representation 
$$
\Ind_{\mu_n(k)\times \pi_1(\und{\Bun}_{T^{\sharp}})}^{\bar\Gamma} (\zeta\boxtimes\chi)\ .
$$
Pick a vector space $\VV$ and an isomorphism of the above representation with $\VV\otimes \Ind_{\bar H}^{\bar\Gamma} \bar\chi$. It yields a decomposition $\pi_! W\iso \VV\otimes\cK_E$ over $\Bunt_{T,\lambda}^{\mu}$ for each $\mu\in\Lambda^{\sharp}$. The Hecke property of $\cK_E$ is obtained from that of $\pi_! W$.
\end{Prf}

\begin{Rem}  In our setting what really matters in the input data $(\theta, n)$ is the bilinear form $\frac{\kappa}{n}\colon \Lambda\otimes\Lambda\to \QQ/\ZZ$ (compare with \cite{FL}, Remark~1). Namely, let $d\ge 1$. Assume $nd$ invertible in $k$.  Let $\bar\zeta\colon \mu_{nd}(k)\to\Qlb^*$ be a character satisfying $\bar\zeta^d=\zeta$. Let $\Bunt_{T, \lambda^d}$ be the gerbe of $nd$-th roots of $\lambda^d$. We have a morphism 
$$
f_d\colon \Bunt_{T,\lambda}\to\Bunt_{T, \lambda^d}
$$ 
sending $(\cB, \cF, \cB^n\iso \lambda_{\cF})$ to $(\cB, \cF, \cB^{nd}\iso\lambda^d_{\cF})$. The functor $f_d^*$ gives rise to an equivalence 
$$
f_d^*\colon \D_{\bar\zeta}(\Bunt_{T, \lambda^d})\iso \D_{\zeta}(\Bunt_{T, \lambda}) \ .
$$
So, $(\theta, n)$ and $(\theta^d, dn)$ give rise essentially to the same problem of the spectral decomposition.
\end{Rem}

\subsubsection{} By Lemma~\ref{Lm_K-torsors_equivariance}, for $\cG_1,\cG_2\in \Bun_K$ we have naturally $\tau_{\cG_1}\otimes\tau_{\cG_2}\iso \tau_{\cG_1\otimes\cG_2}$. Let $\Lambda(\tau)$ denote the biextension of $\Bun_{T^{\sharp}}\times\Bun_K$, whose fibre at $\cG, \cG_1$ is 
$$
\tau_{\cG\otimes\cG_1}\otimes \tau_{\cG}^{-1}\otimes\tau_{\cG_1}^{-1} \ .
$$
Lemma~\ref{Lm_K-torsors_equivariance} yields a trivialization of this biextension. Now by (\cite{P}, Theorem~10.5), there is a line bundle $\tau_H$ on 
$\Bun_{T^{\sharp}, H}$, whose restriction to $'\Bun_{T^{\sharp}}$ is identified with $'\tau$. 

\medskip\noindent
{\bf Question.} Is it true that the line bundle $\tau_H^n\otimes \lambda^{-1}$ is trivial on $\Bun_{T^{\sharp}, H}$? We know already that the restriction of  the central extension (\ref{ext_H^1(X,K)_by_mu_n}) to $H$ is abelian. Does the restriction of  the central extension (\ref{ext_H^1(X,K)_by_mu_n}) to $H$ split? Is it true that the gerbe $\Bunt_{T^{\sharp}, H}\to \Bun_{T^{\sharp}, H}$ is trivial?  

\subsubsection{Example (1)}
\label{section_case_T=Gm}

 Take $\Lambda=\ZZ$, so $T=\Gm$ and $\kappa\colon\Lambda\otimes\Lambda\to\ZZ$ given by $\kappa(x_1,x_2)=2x_1x_2$. There is an object $\theta=(\kappa,\lambda,c)\in \cP^{\theta}(X,\Lambda)$ such that the corresponding line bundle $\lambda$ on $\Bun_T$ has fibre $\det\RG(X, L)\otimes\det\RG(X, L^{-1})\otimes \det\RG(X,\cO)^{-2}$ at $L\in\Bun_1$. We get $\delta_{\lambda}=\kappa$ here. 
Pick a square root of $(\kappa,\lambda,c)$ in $\cP^{\theta}(X,\Lambda)$, it yields as in Section~\ref{subsection_theta_data} a super line bundle $\cL_1$ on $\Bun_1$ with a $\ZZ/2\ZZ$-graded isomorphism $\cL_1^2\iso \lambda$ on $\Bun_T$. Over $\Bun_1^0$ the line bundle $\cL_1$ descends to $\uBun_1^0$. 

 Let $n\ge 1$. Let $e\ge 1$ be the smallest positive integer such that $2e\in n\ZZ$. So, $e=n$ for $n$ odd and $e=n/2$ for $n$ even. We get $\Lambda^{\sharp}=e\ZZ$. Identify $T^{\sharp}$ with $\Gm$, so that $i_X\colon \Bun_1=\Bun_{T^{\sharp}}\to\Bun_T=\Bun_1$ is $L\mapsto L^e$.
 
\begin{Lm} For any $L\in\Bun_1$ there is a canonical $\ZZ/2\ZZ$-graded isomorphism
$$
\frac{\det\RG(X, L^e)\otimes\det\RG(X, L^{-e})}{\det\RG(X, L)^{e^2}\otimes \det\RG(L^{-1})^{e^2}}\iso \det\RG(X,\cO)^{2-2e^2}\ .
$$
\end{Lm}
\begin{Prf}
Let $K(L)$ denotes the LHS of the formula to be proved. Note that $K(L)$ descends to a line bundle on $\uBun_1$. It suffices to check that for $x\in X$ one has $K(L(x))\iso K(L)$ canonically. Since $\det\RG(X, \cdot)$ is multiplicative in exact sequences of coherent sheaves on $X$, this is reduced to showing that 
$$
\frac{\det\RG(X, L^e(ex)/L^e)\otimes L_x^{-e^2}}{\det\RG(X, L^{-e}/L^{-e}(-ex))\otimes (L(x)/L)^{e^2}}
$$
is canonically trivialized. Using the fact that for a line bundle $\cA$ on $X$ one has canonically $\det\RG(X, \cA(ex)/\cA)\iso \cA_x^e\otimes\det\RG(X, \cO(ex)/\cO)$, our claim is reduced to a canonical isomorphism
$\det\RG(X, \cO(ex)/\cO)\iso \det\RG(X, \cO/\cO(-ex))\otimes \Omega_x^{-e^2}$. 
\end{Prf}
 
\medskip

 We may pick (\ref{iso_nth_root_from_lambda}) here with the following properties. The line bundle $\tau$ on $\Bun_{T^{\sharp}}=\Bun_1$ is as follows. If $n$ is odd then $\tau_L=\det\RG(X, L)^n\otimes\det\RG(X, L^{-1})^n\otimes\det\RG(X, \cO)^{-2n}$. If $n$ is even then $\tau=\cL_1^e$ is a super line bundle.  
 
 The group $K(\tau)\subset \uBun_1^0$ is as follows. If $n$ is odd then $K(\tau)=\H^1(X, \mu_{2n})$, and the inclusion $\H^1(X, K)\subset K(\tau)$ is strict. If $n$ is even then $K(\tau)=\H^1(X, \mu_e)$, and $\H^1(X, K)=K(\tau)$. 

\medskip\noindent
\select{Example (2).}
Consider the case when $\kappa(\Lambda\otimes\Lambda)\subset 2\ZZ$. Pick a presentation $\frac{\kappa}{2}=\sum_i b_i(\check{\lambda}_i\otimes\check{\lambda}_i)$ with $\check{\lambda}_i\in \check{\Lambda}$. For $\check{\lambda}\in\check{\Lambda}$ let $R^{\check{\lambda}}$ be the line bundle on $\Bun_T$ with fibre 
$$
\frac{\det\RG(X, \cL^{\check{\lambda}}_{\cF})\otimes \det\RG(X, \cL^{-\check{\lambda}}_{\cF})}{\det\RG(X, \cO)^2}
$$
at $\cF\in\Bun_T$. There is $\theta=(\kappa,\lambda,c)\in\cP^{\theta}(X,\Lambda)$ such that the corresponding line bundle $\lambda$ on $\Bun_T$ is $\lambda=\bigotimes_i (R^{\check{\lambda}_i})^{b_i}$, and we get $\delta_{\lambda}=\kappa$. 

\subsubsection{Corrected scalar products} Let $E, E'$ be 
$\check{T}^{\sharp}$-local systems on $X$. The local system $AE$ on $\Bun_{T^{\sharp}}$ descends to a (defined up to a unique isomorphism) local system on $\und{\Bun}_{T^{\sharp}}$ that we denote by the same letter by abuse of notation. For each $\mu\in \Lambda^{\sharp}$ we may consider the `corrected scalar product' of $AE$ and $AE'$, namely 
\begin{equation}
\label{scalar_product_for_Bun_T_sharp}
\RG(\und{\Bun}_{T^{\sharp}}^{\mu}, AE^*\otimes AE').
\end{equation}
The word `corrected' here refers to the rigidification of $\Bun_T$ along $T$. A similar rigidification (along the center of $\GL_n$) has appeared in the calculation of the scalar product of automorphic sheaves for $\GL_n$ in \cite{L}. 
\begin{Lm} The complex (\ref{scalar_product_for_Bun_T_sharp}) vanishes unless $E\iso E'$. If $E\iso E'$ then for each $\mu\in\Lambda^{\sharp}$ the complex (\ref{scalar_product_for_Bun_T_sharp}) identifies canonically with $\bigoplus_i \wedge^i\H^1(\und{\Bun}_T^0,\Qlb)[-i]$. \QED
\end{Lm}
 
 To express the dependence of the automorphic sheaf $\cK$ of Proposition~\ref{Pp_H_eigensheaves} on $E$, let us write $\cK_E=\cK$. Consider the complex $(\cK_E)^*\otimes\cK_{E'}$ on $\Bunt_{T,\lambda}$. Note that $\mu_n(k)$ acts trivially on this complex, so it descends to a complex (defined up to a unique isomorphism) on $\und{\Bun}_T$, we denote it by the same letter by abuse of notation. The `corrected scalar product' of $\cK_E$ and $\cK_{E'}$ is
\begin{equation}
\label{scalar_product_for_Bun_T}
\RG(\und{\Bun}_T^{\mu}, (\cK_E)^*\otimes\cK_{E'}).
\end{equation}

\begin{Pp} 
\label{Pp_scalar_products}
For each $\mu\in\Lambda^{\sharp}$ one has canonically
\begin{equation}
\label{iso_scalar_products}
\RG(\und{\Bun}_T^{\mu}, (\cK_E)^*\otimes\cK_{E'})\iso 
\RG(\und{\Bun}_{T^{\sharp}}^{\mu}, AE^*\otimes AE').
\end{equation}
\end{Pp}
\begin{Prf}
If $E$ and $E'$ are not isomorphic the both sides vanish, so we assume $E=E'$. Recall that the map $\pi_H\colon \Bunt_{T^{\sharp}, H}\to \Bunt_{T,\lambda}$ defined in Section~\ref{Section_524} is an \'etale Galois covering with Galois group $\H^1(X, K)/H$. Let $AE_H$ and $\bar\chi$ be as in the proof of Proposition~\ref{Pp_H_eigensheaves}. Note that $\pi_H^*\cK_E\iso \bigoplus_{\sigma\in \H^1(X, K)/H} \;\sigma^*AE_H$, so 
$$
(\cK_E)^*\otimes \cK_E\iso \pi_{H !}(\bigoplus_{\sigma\in \H^1(X, K)/H}\; AE_H^*\otimes \sigma^*AE_H)\ .
$$
The local system $\sigma^*AE_H$ corresponds to the character $\bar\chi^{\sigma}$ obtained from $\bar\chi$ by conjugating with any element $\bar\sigma\in \bar\Gamma$ over $\sigma$. 
 Write $\und{\Bun}_{T^{\sharp}, H}$ for the coarse moduli space of $\Bun_{T^{\sharp}, H}$. For each $\sigma$, the local system $AE_H^*\otimes \sigma^*AE_H$ descends with respect to the gerbe $\Bunt_{T^{\sharp}, H}\to \und{\Bun}_{T^{\sharp}, H}$. For $\sigma\in  \H^1(X, K)/H$ we have $\bar\chi=\bar\chi^{\sigma}$ if and only if $\sigma=0$, because $(\cdot,\cdot)_c$ is non-degenerate. So, only $\sigma=0$ contributes nontrivially to the LHS of (\ref{iso_scalar_products}), which idenitifies with $\RG(\und{\Bun}^{\mu}_{T^{\sharp}, H}, \Qlb)$.
\end{Prf}

\subsubsection{About Question~\ref{Con_one}}
\label{Section_528} If the answer to Question~\ref{Con_one} is positive then for any $\mu\in\Lambda^{\sharp}$ there is $M\in \D_{\zeta}(\Bunt_{T,\lambda}^{\mu})$ such that for any $\check{T}^{\sharp}$-local system $E$ on $X$ the complex $\R\Hom(\cK_E, M)$ is placed in cohomological degree zero and is of dimension one. 

 Here is a model situation for $T=\Gm$ and $g=1$ showing that one should expect the negative answer to Question~\ref{Con_one}. Let $Y, X$ be elliptic curves and $i\colon Y\to X$ be an isogeny with kernel $K$, which is reduced of order $e^2$.  Assume given a central extension $1\to \mu_n(k)\to \Gamma\to K\to 1$ such that the corresponding commutator pairing yields an isomorphism $K\iso \Hom(K, \mu_n(k))$. Let $\Gamma$ act on $Y$ via its quotient $K$, write 
$\tilde X$ for the stack quotient of $Y$ by $\Gamma$. Pick an injective character $\zeta\colon \mu_n(k)\to\Qlb^*$. Let $\D_{\zeta}(\tilde X)$ be the bounded derived category of $\Qlb$-sheaves on $\tilde X$ on which $\mu_n(k)$ acts by $\zeta$.
Each irreducible local system $E\in \D_{\zeta}(\tilde X)$ is of rank $e$. 

\begin{Pp} Let $E$ be an irreducible local system in $\D_{\zeta}(\tilde X)$. If $M\in \D_{\zeta}(\tilde X)$ then $e$ divides $\chi(\Spec k, \R\Hom(E, M))$.
\end{Pp}
\begin{Prf}
Let $h\colon Y\to \tilde X$ be the quotient map. By (\cite{FGV}, Theorem~6.7), we get $\chi(Y, h^*(E^*\otimes M))=e\chi(Y, h^*M)$, as we may replace $h^*E$ by $\Qlb^e$. The complex $h_*\Qlb$ decomposes by the characters of $\mu_n(k)$, and only the trivial character of $\mu_n(k)$ contributes to $\chi(\tilde X, (E^*\otimes M)\otimes h_*\Qlb)=e^2\chi(\tilde X, (E^*\otimes M))$. So, $\chi(Y, h^*M)=e\chi(\tilde X, E^*\otimes M)$. 

 Let us show that if $F\in \D_{\zeta}(\tilde X)$ then $\chi(Y, h^*F)$ is divisible by $e^2$. Indeed, $\chi(Y, h^*F)$ depends only on the image of $F$ in the Grothendieck group of $\D_{\zeta}(\tilde X)$. The latter is generated by the irreducible perverse sheaves. Let $F\in \D_{\zeta}(\tilde X)$ be an irreducible perverse sheaf. If $F$ is supported on the preimage of a point in $X$ then this is clear. Assume now $F$ is supported generically. For each point $x\in X$, where $F$ is not a local system, $i^{-1}(x)$ consists of $e^2$ elements. By the Ogg-Shafarevich formula for the Euler characteristic, the local contributions in $\chi(Y, h^*F)$ for all points in $i^{-1}(x)$ are the same. So, $\chi(Y, h^*F)$ is divisible by $e^2$.
\end{Prf}

\medskip\noindent
\select{Acknowledgements.} When preparing this paper the author learned that 
similar results in the context of $\cD$-modules were obtained several years ago by D. Gaitsgory (unpublished) in the guise of quantum geometric Langlands correspondence for the torus, we are grateful to him for useful discussions.  We also thank the organizers of the conference `Automorphic forms and harmonic analysis on covering groups' (The American Institute of Mathematics, Palo Alto, 2013), which stimulated our work. We are also grateful to V. Lafforgue and 
M.~Finkelberg for many fruitful discussions. The author was supported by ANR-13-BS01-0001-01 grant.

\end{document}